\documentclass[preprint,12pt]{elsarticle}

\usepackage{amssymb}
\usepackage{amsthm}
\usepackage{framed} 
\usepackage{amsmath,color}
\usepackage{mathrsfs}
\usepackage{graphicx}
\usepackage{epstopdf}
\usepackage{float}
\usepackage{caption}
\usepackage{subcaption}
\usepackage{bm}
\usepackage{bbm}
\usepackage{mathrsfs}
\usepackage{cleveref}
\usepackage{soul}
\usepackage{accents}
\usepackage{color,soul} 
\usepackage{color} 
\usepackage{nomencl} 
\makenomenclature 
\biboptions{sort&compress}
\soulregister\citep7 
\soulregister\citet7 
\soulregister\citealp7 
\newsavebox{\measurebox} 
\usepackage{titlesec} 
\usepackage{tabu} 
\usepackage{longtable}
\usepackage[euler]{textgreek} 

\journal{ }

\makeatletter
\def\@author#1{\g@addto@macro\elsauthors{\normalsize%
    \def\baselinestretch{1}%
    \upshape\authorsep#1\unskip\textsuperscript{%
      \ifx\@fnmark\@empty\else\unskip\sep\@fnmark\let\sep=,\fi
      \ifx\@corref\@empty\else\unskip\sep\@corref\let\sep=,\fi
      }%
    \def\authorsep{\unskip,\space}%
    \global\let\@fnmark\@empty
    \global\let\@corref\@empty  
    \global\let\sep\@empty}%
    \@eadauthor={#1}
}
\makeatother

\titleformat{\paragraph}
{\normalfont\normalsize\itshape}{\theparagraph}{1em}{}
\titlespacing*{\paragraph}
{0pt}{3.25ex plus 1ex minus .2ex}{1.5ex plus .2ex}

\begin{document}

\begin{frontmatter}



\title{A phase field formulation for hydrogen assisted cracking}


\author{Emilio Mart\'{\i}nez-Pa\~neda\corref{cor1}\fnref{Cam}}
\ead{mail@empaneda.com}

\author{Alireza Golahmar\fnref{DTU}}

\author{Christian F. Niordson\fnref{DTU}}

\address[Cam]{Department of Engineering, University of Cambridge, CB2 1PZ Cambridge, UK}

\address[DTU]{Department of Mechanical Engineering, Technical University of Denmark, DK-2800 Kgs. Lyngby, Denmark}

\cortext[cor1]{Corresponding author.}

\begin{abstract}
We present a phase field modeling framework for hydrogen assisted cracking. The model builds upon a coupled mechanical and hydrogen diffusion response, driven by chemical potential gradients, and a hydrogen-dependent fracture energy degradation law grounded on first principles calculations. The coupled problem is solved in an implicit time integration scheme, where displacements, phase field order parameter and hydrogen concentration are the primary variables. We show that phase field formulations for fracture are particularly suitable to capture material degradation due to hydrogen. Specifically, we model (i) unstable crack growth in the presence of hydrogen, (ii) failure stress sensitivity to hydrogen content in notched specimens, (iii) cracking thresholds under constant load, (iv) internal hydrogen assisted fracture in cracked specimens, and (v) complex crack paths arising from corrosion pits. Computations reveal a good agreement with experiments, highlighting the predictive capabilities of the present scheme. The work could have important implications for the prediction and prevention of catastrophic failures in corrosive environments. The finite element code developed can be downloaded from www.empaneda.com/codes
\end{abstract}


\end{frontmatter}



\section{Introduction}
\label{Sec:Introduction}

Hydrogen has been known for over a hundred years to cause catastrophic failure in metallic materials \cite{Johnson1875}. Hydrogen atoms enter the material, migrate through the crystal lattice and induce fracture, with cracking being observed in modern steels at one-tenth of the expected fracture toughness. Hydrogen embrittlement is a pervasive problem in the manufacturing, transport, and energy sectors and could compromise the key role of hydrogen as next-generation energy resource. There is a strong need to understand and model how hydrogen degrades the ductility and the fracture resistance of metals.\\

Modeling efforts in the hydrogen embrittlement community have been largely restricted to discrete numerical methods. Crack initiation thresholds and growth rates have been computed through dislocation-based methods \cite{AM2016,Burns2016}, weakest-link approaches \cite{Novak2010,Ayas2014}, and, eminently, cohesive zone models \cite{Serebrinsky2004,Scheider2008,Moriconi2014,EFM2017,Yu2017}, which can easily incorporate the cohesive strength reduction with hydrogen coverage. These methods are however limited when dealing with the complex conditions of practical applications and many important challenges (mixed mode, interacting cracks, complex distribution of initiation sites, etc.) remain unaddressed.\\

The phase field method for fracture has emerged as a promising variational approach to overcome the limitations of discontinuity-based methods. The method has received increasing attention since the earlier work by Francfort and Marigo \cite{Francfort1998} and important efforts have been devoted to improve solution schemes \cite{Miehe2010,Gerasimov2016}, discretization strategies \cite{Borden2014,Schillinger2015} and implementation aspects \cite{Bourdin2000,Reinoso2017}. The model builds upon Griffith's thermodynamics framework and has not only been successfully applied to model brittle fracture but also to ductile damage \cite{Borden2016,Miehe2016a}, hydraulic fracturing \cite{Mikelic2015,Wilson2016}, composites delamination \cite{Reinoso2017a,Carollo2017}, and stress corrosion cracking \cite{Mai2016,Mai2017}, to name a few. The phase field method is particularly well suited to capture the fracture energy dependence on hydrogen concentration, and will allow addressing the complex cracking conditions intrinsic to hydrogen embrittlement in industrial components.\\

In this work, we develop a novel variational framework for hydrogen assisted cracking. Hydrogen transport towards the fracture process zone and subsequent cracking are investigated by means of a coupled mechanical-diffusion-phase field finite element framework. The numerical scheme is grounded on: (i) a phase field description of fracture, (ii) an extension of Fick's law for mass diffusion, and (iii) a quantum-based relation between the hydrogen content and the fracture surface energy. We demonstrate the potential of the proposed modeling framework by means of representative case studies.nThe remainder of the paper is organized as follows. Section \ref{Sec:PhaseField} describes the hydrogen-dependent phase field formulation. Section \ref{Sec:GeneralFramework} provides details of the stress-assisted impurity diffusion scheme and its coupling with the phase field model. In Section \ref{Sec:Results} model predictions are benchmarked against relevant experimental measurements. Finally, concluding remarks are given in Section \ref{Sec:Concluding remarks}.

\section{A phase field formulation for hydrogen embrittlement}
\label{Sec:PhaseField}
\subsection{Phase field approximation of fracture}
\label{Sec:PhaseFieldFracture}

\begin{figure}[H] 
    \centering
    \includegraphics[scale=0.9]{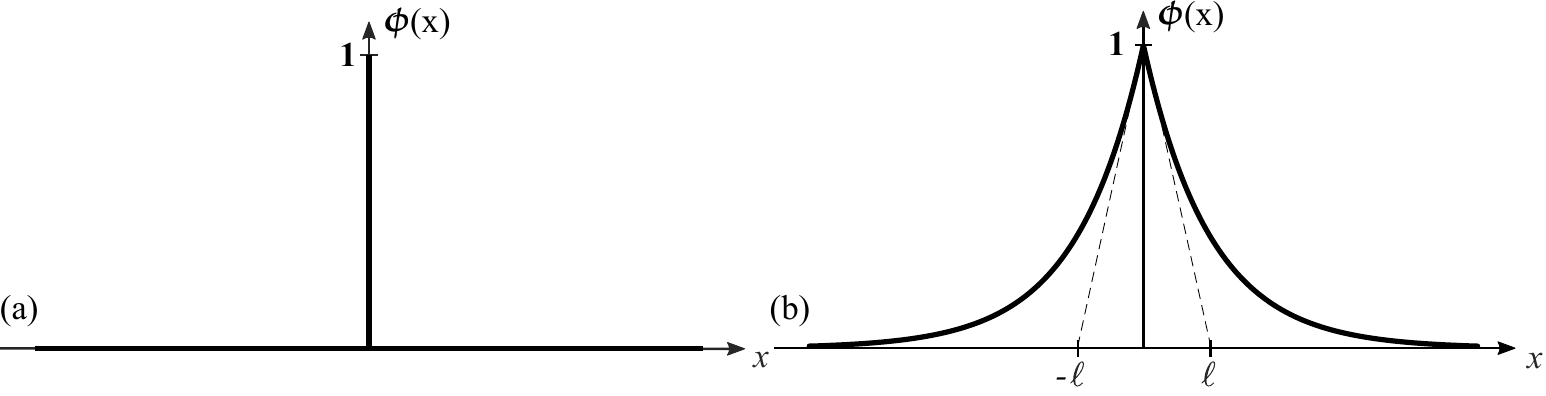}
    \caption{Sharp (a) and diffusive (b) crack topologies.}
    \label{fig:SharpSmooth}
\end{figure}

In a one-dimensional setting, the topology of a sharp crack (see Fig. \ref{fig:SharpSmooth}a) can be described by an auxiliary field variable \(\phi(x) \in[0,1]\) with
\begin{equation} \label{eSharp}
    \phi(x) =
    \begin{cases}
      1  & \quad \text{if } x=0\\
      0  & \quad \text{if } x \neq 0\\
    \end{cases}
    \centering
\end{equation}
which is referred to as the crack phase field order parameter, with \(\phi=0\) and \(\phi=1\) respectively denoting the intact and fully broken states of the material. The non-smooth crack phase field (\ref{eSharp}) can be approximated by the exponential function,
\begin{equation}
    \phi(x)= \text{e}^{-\tfrac{|x|}{\ell}}
    \centering
    \label{eDiffusive}
\end{equation}
\noindent representing a regularized or diffusive crack topology, as illustrated in Fig. \ref{fig:SharpSmooth}b. The length-scale parameter \(\ell\) determines the width of the smearing function, approaching the discrete crack topology as \(\ell \rightarrow 0\). 

\begin{figure}[H] 
    \centering
    \includegraphics[scale=1]{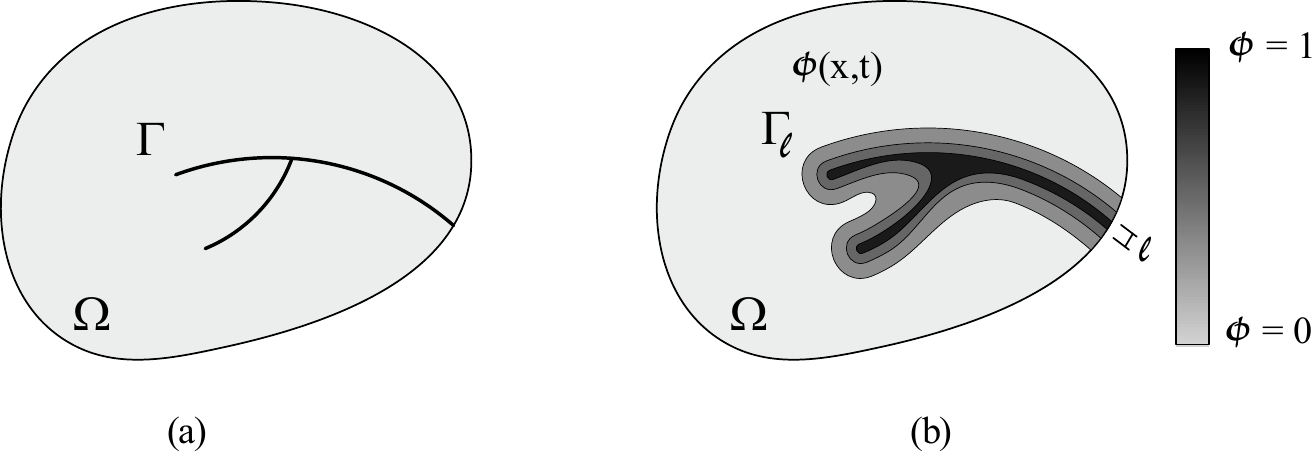}
    \caption{Schematic representation of a solid body with (a) internal discontinuity boundaries, and (b) a phase field approximation of the discrete discontinuities. Adapted from \cite{Borden2012}.}
    \label{fig:PhaseFieldFigIntro}
\end{figure}

Consider a discrete internal discontinuity $\Gamma$ in a solid body $\Omega$ (Fig. \ref{fig:PhaseFieldFigIntro}a), a regularized crack functional $\Gamma_\ell$ can be defined (Fig. \ref{fig:PhaseFieldFigIntro}b)
\begin{equation}
    \Gamma_{\ell}(\phi)= \int_{\Omega} \left( \dfrac{1}{2\ell}\phi^{2} + \dfrac{\ell}{2}{|\nabla \phi|}^{2} \right) \, \mathrm{d}V = \int_{\Omega} \gamma_\ell(\phi,\nabla \phi) \, \mathrm{d}V
    \centering
\end{equation}
\noindent that \textGamma -converges to the functional of the discrete crack surface for a vanishing length scale parameter $\ell \to 0$. Here, \(\gamma_\ell(\phi,\nabla \phi)\) is the crack surface density function. Consequently, the phase field method for fracture circumvents numerical complications inherent to tracking the evolving discontinuity boundary $\Gamma$, enabling to robustly model interactions and branching of cracks of arbitrary topological complexity \cite{Borden2012}. The fracture energy due to the formation of a crack is then approximated as,
\begin{equation}\label{Eq:Gcapproximation}
    \int_\Gamma G_c \left( \theta \right) \mathrm{d}S \approx \int_{\Omega} G_c \left( \theta \right) \left( \dfrac{1}{2\ell}\phi^{2} + \dfrac{\ell}{2}{|\nabla \phi|}^{2} \right) \, \mathrm{d}V
\end{equation}
\noindent where $G_c$ is the critical Griffith-type energy release rate, which is dependent on the hydrogen coverage $\theta$.

The phase field fracture method resembles traditional continuum damage mechanics models, where the scalar damage field may be interpreted as the phase field $\phi$. Accordingly, the loss of stiffness associated with mechanical degradation of the material is characterized as a function of $\phi$. This is done by means of the so-called degradation function $g \left( \phi \right)$, which relates the stored bulk energy per unit volume to the strain energy density of the undamaged solid,
\begin{equation}\label{Eq:psi}
\psi=g \left( \phi \right) \psi_0 = \left[ \left( 1 - \phi \right)^2 + k \right] \psi_0
\end{equation} 
Here, $k$ is a parameter chosen to be as small as possible to keep the system of equations well-conditioned; a value of $k=1 \times 10^{-7}$ is adopted throughout this work. Integrating over the volume of the body and adding the phase field free energy contribution (\ref{Eq:Gcapproximation}) one reaches the total potential energy functional of the deformation-fracture problem,
\begin{equation}\label{Eq:PsiDefPhaseField}
    \Psi=\int_{\Omega} \left\{ \left[(1-\phi)^{2}+k\right] \psi_0 + G_c \left( \theta \right) \left(\dfrac{1}{2\ell}\phi^{2} + \dfrac{\ell}{2}{|\nabla \phi|}^{2}\right) \right\} \, \mathrm{d}V
    \centering
\end{equation}
The variation of (\ref{Eq:PsiDefPhaseField}) with respect to $\delta \phi$ renders the weak form of the phase field contribution,
\begin{equation}\label{Eq:weak}
  \int_{\Omega} \left[ -2(1-\phi)\delta \phi \, \psi_{0} +
        G_c \left( \theta \right) \left( \dfrac{1}{\ell}\phi \delta \phi
        + \ell\nabla \phi \cdot \nabla \delta \phi \right) \right]  \, \mathrm{d}V = 0
\end{equation}
and the phase field equilibrium equation in $\Omega$,
\begin{equation} \label{Eq:strong}
G_c \left( \theta \right) \left(\dfrac{1}{\ell} \phi - \ell \Delta \phi \right) - 2(1-\phi) \, \psi_{0} = 0
\end{equation}
follows immediately upon making use of the product rule and Gauss' divergence theorem.

\subsection{Hydrogen-dependent surface energy degradation}
\label{Sec:Henergydegradation}

\begin{figure}[H]
\centering
\includegraphics[scale=1]{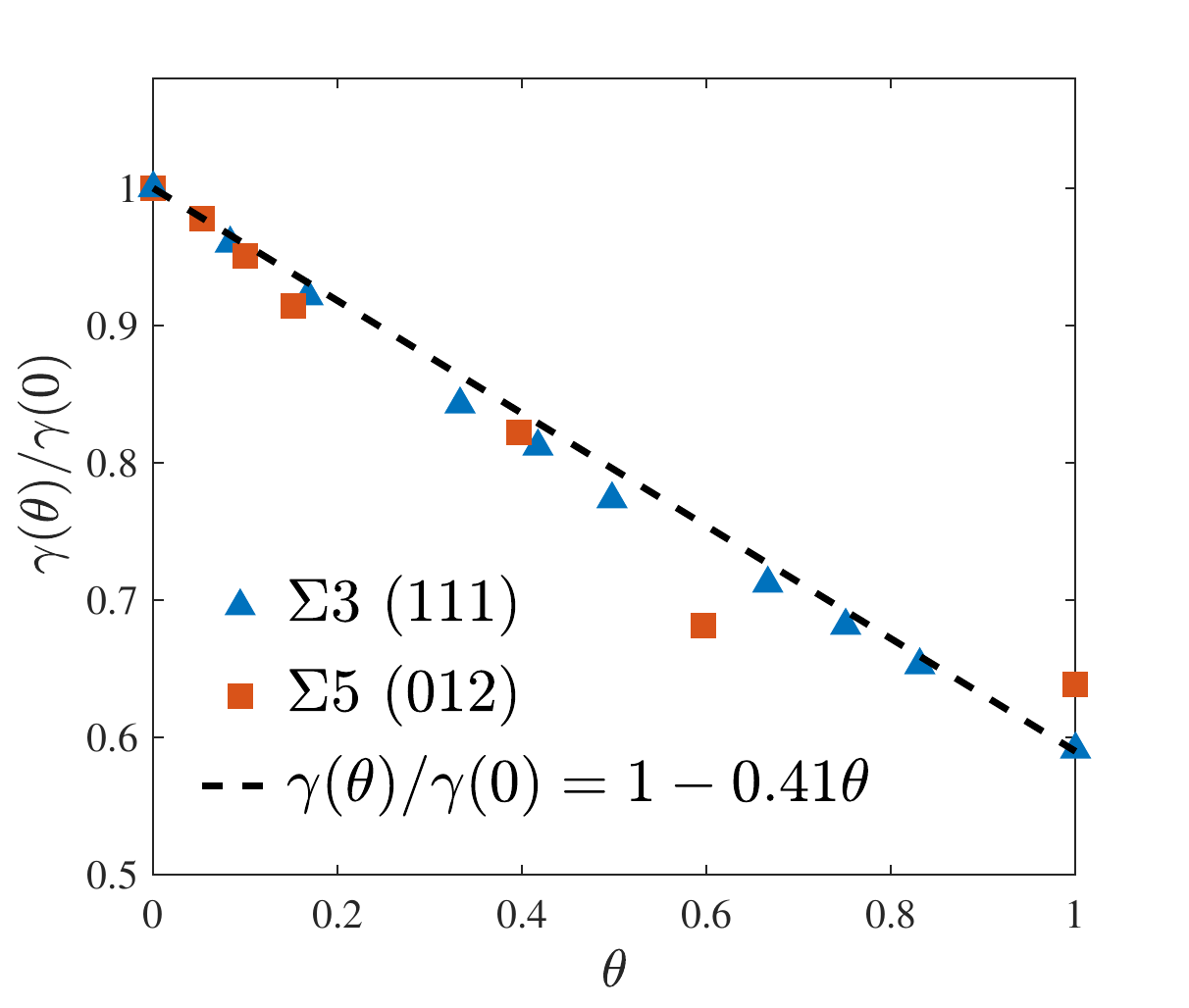}
\caption{Effect of hydrogen on the surface energy of nickel. Linear fit to the DFT calculations by Alvaro et al. \cite{Alvaro2015} for $\Sigma3$ and $\Sigma5$ grain boundaries.}
\label{fig:DFT}
\end{figure}

Hydrogen embrittles metallic materials by lowering the bond energy between metal atoms, which translates into a reduction of the fracture resistance. A number of authors have employed Density Functional Theory (DFT) to investigate the decohesion of fracture surfaces with varying hydrogen coverage (see, e.g., \cite{Jiang2004a,Alvaro2015,Tehranchi2017a} and references therein). For example, Alvaro et al. \cite{Alvaro2015} computed the change in ideal fracture energy in the presence of hydrogen atoms at $\Sigma3$ and $\Sigma5$ grain boundaries in nickel. Their results, in terms of normalized surface energy $\gamma \left( \theta \right) / \gamma \left( 0 \right)$ versus hydrogen coverage $\theta$, are shown in Fig. \ref{fig:DFT}, along with a linear fit to the data.

Equivalently, one can define the critical energy release rate dependence on the hydrogen coverage as,
\begin{equation}
    \frac{G_c \left( \theta \right)}{G_c \left( 0\right)}=1 - \chi \theta
\end{equation}

\noindent where $G_c\left( 0\right)$ is the critical energy release rate in the absence of hydrogen and $\chi$ is the damage coefficient that weights the hydrogen-lowering of the fracture energy. $\chi$ can be estimated for other materials by fitting DFT data from the literature. The damage coefficient for iron and aluminum is obtained from the work by Jiang and Carter \cite{Jiang2004a} - see Table \ref{Tab:DamageCoeff}.

\begin{table}[H]
\centering
\caption{Weighting factor estimation from first principles quantum mechanics.}
\label{Tab:DamageCoeff}
   {\tabulinesep=1.2mm
   \begin{tabu} {ccc}
       \hline
 Material & Damage coefficient $\chi$ & DFT analysis \\ \hline
 Iron   & 0.89 & Jiang and Carter (2004) \cite{Jiang2004a}\\
 Nickel & 0.41 & Alvaro et al. (2015) \cite{Alvaro2015} \\
 Aluminum & 0.67 & Jiang and Carter (2004) \cite{Jiang2004a}\\\hline
   \end{tabu}}
\end{table}

The effect of hydrogen on the fracture resistance can be illustrated by examining the analytical homogeneous solution of a one-dimensional quasi-static problem. Under these circumstances, the Cauchy stress is given by $\sigma=g(\phi)E \epsilon$, where $E$ denotes Young's modulus and $\epsilon$ the strain. Given a strain energy density $\psi_0=E \epsilon^2/2$, one can readily obtain the homogeneous phase field from the strong form (\ref{Eq:strong}),
\begin{equation}
    \phi=\frac{E \epsilon^2 \ell}{G_c \left( \theta \right)+E \epsilon^2 \ell}
\end{equation}

\noindent and substituting into the constitutive equation renders the characteristic relation between the homogeneous strain and the homogeneous stress,
\begin{equation}
    \sigma=\left( \frac{G_c \left( \theta \right)}{G_c\left( \theta \right) + E \varepsilon^2 \ell} \right)^2 E \varepsilon
\end{equation}

The homogeneous solution for the stress reaches a maximum at a critical stress quantity,
\begin{equation}\label{Eq:SigmaC}
    \sigma_c= \sqrt{\frac{27 E G_c \left( \theta \right)}{256 \ell}}
\end{equation}

\noindent with the strain counterpart given by,
\begin{equation}
    \epsilon_c = \sqrt{\frac{G_c \left( \theta \right)}{3 \ell E}}
\end{equation}

\begin{figure}[H]
\centering
\includegraphics[scale=1]{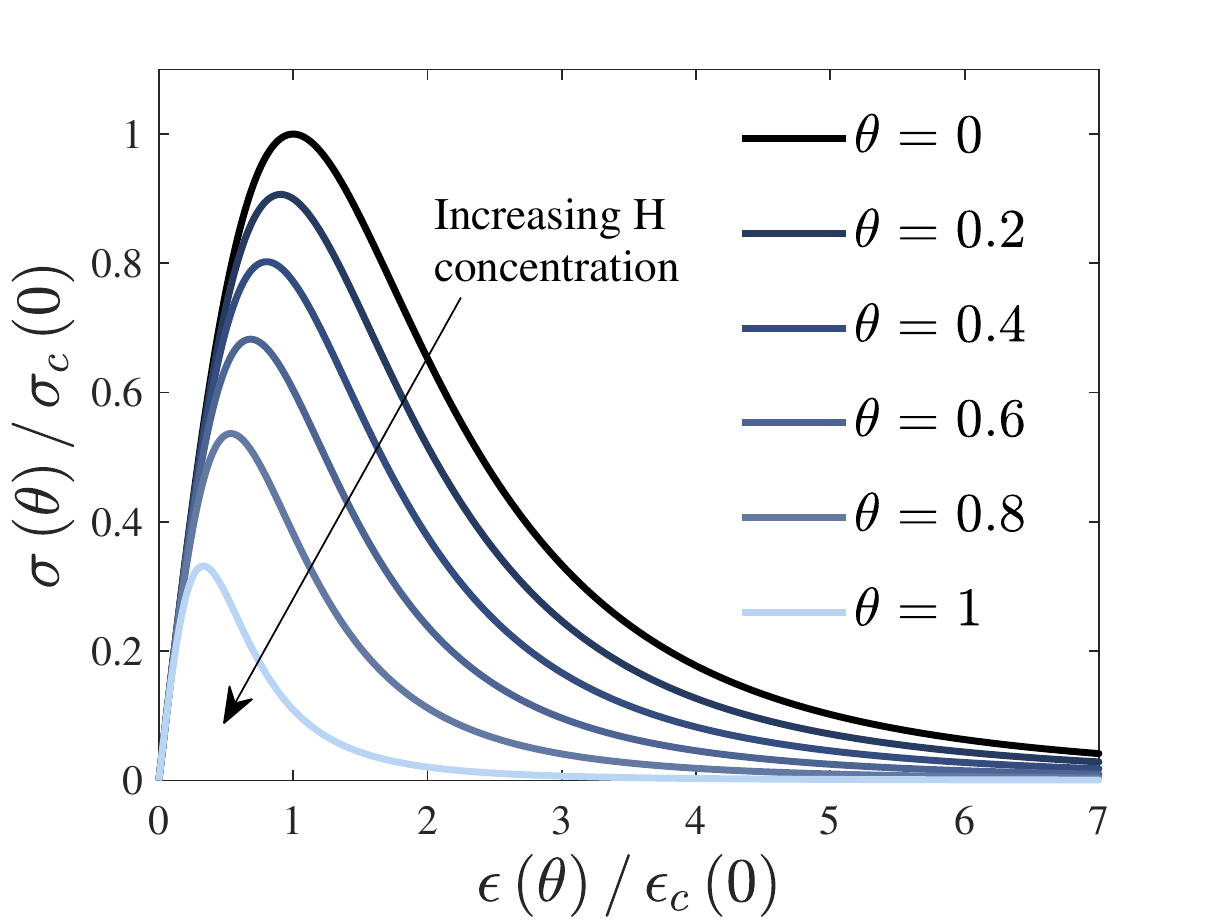}
\caption{Effect of hydrogen coverage $\theta$ on the damage constitutive law for iron-based materials. The stress-strain response is normalized by the corresponding $\theta=0$ quantities.}
\label{fig:Henergy}
\end{figure}

Thus, $\ell$ may be interpreted as a material parameter, which governs the magnitude of the critical stress at which damage initiates. This necessarily implies that the value of $G_c$ is intrinsically related to the choice of $\ell$, being equal to the Griffith's critical energy release rate for $\ell \to 0 $. From (\ref{Eq:SigmaC}) an analogy with cohesive zone formulations can be established. With this in mind, it is interesting to note that the characteristic element size $h$ in cohesive zone analyses is typically chosen to be
\begin{equation}
    h < \frac{\pi}{160} \frac{E G_c}{\sigma_c^2}
\end{equation}

\noindent so as to resolve the fracture process zone with at least 20 elements \cite{Rice1968a,Moes2002}. This implies that $h$ should be at least 5.4 times smaller than $\ell$ in the phase field problem. The constitutive relation of the homogeneous phase field problem is shown in Fig. \ref{fig:Henergy} as a function of hydrogen coverage.

As shown in Fig. \ref{fig:Henergy} for the case of iron, hydrogen significantly degrades the strength and the fracture resistance. The present modeling framework quantitatively accounts for the sensitivity of the surface energy to hydrogen coverage, as observed in first principles calculations; other mechanisms of hydrogen damage can also be incorporated.\\

Finally, we use the Langmuir-McLean isotherm to compute the surface coverage from the bulk hydrogen concentration $C$,
\begin{equation}\label{Eq:Langmuir-McLean}
    \theta=\frac{C}{C+\exp \left( \frac{-\Delta g_b^0}{RT} \right)}
\end{equation}

\noindent with $C$ given in units of impurity mole fraction. Here, $R$ is the universal gas constant, $T$ the temperature and $\Delta g_b^0$ is the Gibbs free energy difference between the decohering interface and the surrounding material. As in \cite{Serebrinsky2004}, a value of 30 kJ/mol is assigned to $\Delta g_b^0$ based on the spectrum of experimental data available for the trapping energy at grain boundaries. Thus, the present formulation accounts for the effect of microstructural traps on cracking and can incorporate the influence on mass transport through an effective diffusion coefficient.

\section{Stress-assisted hydrogen diffusion coupled with phase field fracture}
\label{Sec:GeneralFramework}

\begin{figure}[H] 
    \centering
    \includegraphics[scale=1]{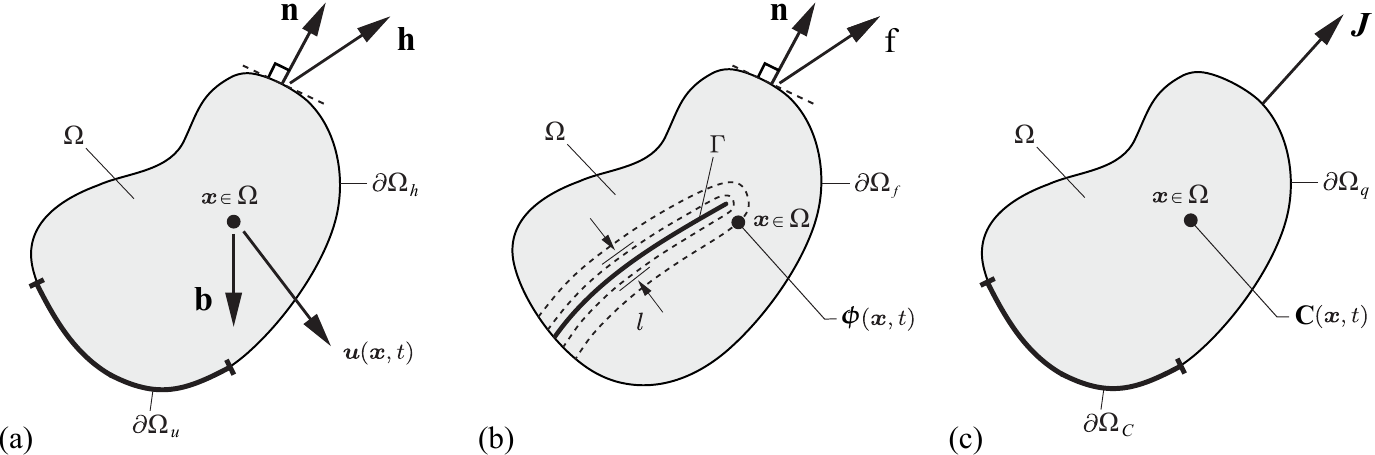}
    \caption{Schematic representation of the three-field boundary value problem: (a) deformation, (b) phase field, and (c) mass transport.}
    \label{fig:BC}
\end{figure}

The formulation presented in this section refers to the response of a solid body $\Omega$ with external surface $\partial \Omega$ of outward normal $\bm{n}$ - see Fig. \ref{fig:BC}. With respect to the displacement field $\bm{u}$, the outer surface of the body is decomposed into a part $\partial \Omega_u$, where the displacement is prescribed by Dirichlet-type boundary conditions, and a part $\partial \Omega_h$, where the traction $\bm{h}$ is prescribed by Neumann-type boundary conditions (see Fig. \ref{fig:BC}a). A body force field per unit volume $\mathbf{b}$ can also be prescribed. The fracture phase field $\phi$ is driven by the displacement field of the solid. A Dirichlet-type boundary condition can be prescribed at $\Gamma$, a given crack surface inside the solid (see Fig. \ref{fig:BC}b). Additionally, a phase field fracture microtraction $f$ can be prescribed on $\partial \Omega_f$. With respect to the hydrogen concentration $C$, the external surface consists of two parts (see Fig. \ref{fig:BC}c): $\partial \Omega_q$, where the hydrogen flux $\bm{J}$ is known (Neumann-type boundary conditions), and $\partial \Omega_C$, where the hydrogen concentration is prescribed (Dirichlet-type boundary conditions). Accordingly, a concentration flux entering the body across $\partial \Omega_q$ can be defined as $q=\bm{J} \cdot \bm{n}$.

\subsection{Governing balance equations}
\label{Sec:PVW}

The principle of virtual work, in the absence of body forces, reads
\begin{align}\label{Eq:PVW}
&\int_\Omega \bigg( \bm{\sigma} : \nabla \delta \bm{u} + \omega \, \delta \phi + \bm{\zeta} \cdot \nabla \delta \phi - \frac{dC}{dt} \, \delta \mu \\ \nonumber
&+ \bm{J} \cdot \nabla \delta \mu \bigg) \, \mathrm{d}V = \int_{\partial \Omega} \left( \bm{h} \cdot \delta \bm{u} + f \, \delta \phi + q \, \delta \mu \right) \, \mathrm{d}S
\end{align}
\noindent where $\mu$ is the chemical potential that drives the flux of hydrogen, $\bm{\sigma}$ is the symmetric Cauchy tensor, and $\omega$ and $\bm{\zeta}$ are microstress quantities work conjugate to the phase field $\phi$ and the phase field gradient $\nabla \phi$, respectively. Applying the product rule and Gauss' divergence theorem to the internal virtual work - left hand side of Eq. (\ref{Eq:PVW}) - renders,
\begin{align}\label{Eq:IVW}
\delta \mathcal{W}_i = &\int_{\partial \Omega} \left( \bm{\sigma} \, \bm{n} \cdot \delta \bm{u} + \bm{\zeta} \cdot \bm{n} \, \delta \phi + \bm{J} \cdot \bm{n} \, \delta \mu  \right) \, \mathrm{d}S - \int_\Omega \bigg( \left( \nabla \cdot \bm{\sigma} \right) \cdot \delta \bm{u}  \\ \nonumber
&+ \left( \nabla \cdot \bm{\zeta} - \omega \right) \delta \phi + \left( \nabla \cdot \bm{J} + \frac{dC}{dt} \right) \delta \mu \bigg) \,\mathrm{d}V 
\end{align}
Since the second integral on the right-hand side of Eq. (\ref{Eq:IVW}) should vanish for arbitrary variations, a set of three equilibrium equations can be obtained,
\begin{equation} \label{strong}
    \begin{split}
        \nabla \cdot \bm{\sigma} &= \mathbf{0} \hspace{5mm} \text{on} \hspace{3mm} \Omega\\[1mm]
        \nabla \cdot \bm{\zeta} - \omega &= 0  \hspace{5mm} \text{on} \hspace{3mm} \Omega\\[1mm]
        \frac{dC}{dt} + \nabla \cdot \bm{J} &= 0  \hspace{5mm} \text{on} \hspace{3mm} \Omega
    \end{split}
    \centering
\end{equation}
The first integral on the right hand side of (\ref{Eq:IVW}) may be identified as part of the external virtual work. Thus, by accounting for the right hand side of (\ref{Eq:PVW}), the corresponding set of boundary conditions is obtained,
\begin{equation} \label{Eq:strongBCs}
    \begin{split}
        \bm{h} = \bm{\sigma} \, \bm{n}  &\hspace{5mm} \text{on} \hspace{3mm} \partial \Omega_h\\[1mm]
        f= \bm{\zeta} \cdot \bm{n}  &\hspace{5mm} \text{on} \hspace{3mm} \partial \Omega_f\\[1mm]
        q= \bm{J} \cdot \bm{n}  &\hspace{5mm} \text{on} \hspace{3mm} \partial \Omega_q 
    \end{split}
    \centering
\end{equation}
\subsection{Energy imbalance}
\label{Sec:Thermodynamics}

The first two laws of thermodynamics within a dynamical process of specific internal energy $\mathscr{E}$ and specific entropy $\Lambda$ read
\begin{equation}
\frac{d}{dt} \int_\Omega \mathscr{E} \, \text{d}V = W_{e}\left(\Omega \right) - \int_{\partial \Omega} \bm{Q} \cdot \bm{n} \, \text{d}S + \int_\Omega Q \, \text{d}V
\end{equation}
\begin{equation}
\frac{d}{dt} \int_\Omega \Lambda \, \text{d}V \geq \int_{\partial \Omega} \frac{\bm{Q}}{\Theta} \cdot \bm{n} \, \text{d}S + \int_{\Omega} \frac{Q}{\Theta} \, \text{d}V
\end{equation}

\noindent where $W_{e}$ is the external work, $\bm{Q}$ is the heat flux, $Q$ is the heat absorption and $\Theta$ the absolute temperature. Within an isothermal process ($\Theta=\Theta_0$) the thermodynamic laws can be readily combined by considering the free energy $\Psi=\mathscr{E} - \Theta_0 \Lambda$, such that
\begin{equation}\label{Eq:ClausiusDuhem}
\frac{d}{dt} \int_\Omega \Psi \, \text{d}V \leq W_{e} \left(\Omega \right)
\end{equation}

\noindent which is generally referred to as Clausius-Duhem inequality or the principle of dissipation \cite{Gurtin2010}. Given (\ref{Eq:PVW}), the local free-energy inequality takes the form
\begin{equation}\label{Eq:ClausiusDuhem}
\frac{d}{dt} \int_\Omega \Psi \, \text{d}V \leq \int_{\partial \Omega_h} \bm{h} \cdot \dot{\bm{u}} \, \text{d}S + \int_{\partial \Omega_f} f \,  \dot{\phi} \, \text{d}S + \int_{\partial \Omega_q} q \, \dot{\mu} \, \text{d}S 
\end{equation}

Employing the divergence theorem and considering (\ref{strong})-(\ref{Eq:strongBCs}) one reaches the local dissipation postulate,
\begin{equation}\label{Eq:Diss}
\bm{\sigma} : \nabla \dot{\bm{u}} + \omega \dot{\phi} + \bm{\zeta} \cdot \nabla \dot{\phi} + \mu \dot{C} + \bm{J} \cdot \nabla \mu  - \dot{\Psi} \leq 0
\end{equation}

\subsection{Constitutive theory}
\label{Sec:ConstitutiveTheory}

The constitutive prescriptions are constructed, in a thermodynamically consistent manner, by ensuring that the dissipation condition (\ref{Eq:Diss}) is satisfied. The free energy is defined as,
\begin{align}\label{Eq:FreeEnergy}
\Psi \left( \bm{\epsilon}, \, \phi, \, \nabla \phi, \, C \right) =& \underbrace{ \left( 1- \phi \right)^2 \psi_0 \left( \bm{\epsilon} \right) - K \bar{V}_H \left( C - C^0 \right) \text{tr} \, \bm{\epsilon}}_{\textstyle \Psi^b} \nonumber  \\[5pt]
& + \underbrace{ G_c \left( \theta \right) \left( \frac{1}{2 \ell} \phi^2 + \frac{\ell}{2} |\nabla \phi|^2 \right)}_{\textstyle \Psi^s} \nonumber  \\[5pt] 
& + \underbrace{\mu^0 C + R T N \left( \theta_L \ln \theta_L + \left( 1 - \theta_L \right) \ln \left( 1 - \theta_L \right) \right)}_{\textstyle \Psi^c}
\end{align}
\noindent where $\Psi^b \left( \bm{\epsilon}, \, \phi, \, C \right)$, $ \Psi^s \left(\phi, \, \nabla \phi, \,  C  \right)$ and $\Psi^c \left( C \right)$ respectively denote the chemo-elastic energy stored in the bulk \cite{DiLeo2013}, the crack surface energy, and the chemical free energy. Here $\bm{\epsilon}$ is the strain tensor, $K$ is the bulk modulus, $\bar{V}_H$ the partial molar volume of hydrogen in solid solution, $C^0$ the reference lattice hydrogen concentration, $\mu^0$ the reference chemical potential, $\theta_L$ the occupancy of lattice sites and $N$ the number of lattice sites. We subsequently develop a constitutive theory from the free energy that couples deformation, fracture and hydrogen diffusion.

\subsubsection{Mechanical deformation}

We define $\psi_0$ as the elastic strain energy density for the undamaged solid,
\begin{equation}
    \psi_{0}(\bm{\epsilon})=\dfrac{1}{2}\bm{\epsilon}^{T} : \mathbf{C_{0}} : \bm{\epsilon}
    \centering
\end{equation}
where \(\textbf{C}_{0}\) is the linear elastic stiffness matrix. A linear elastic description of the material response is considered appropriate to characterize the brittle behavior of metals in the presence of hydrogen. Implicitly multi-scale plasticity formulations (e.g., strain gradient plasticity) have shown that crack tip stresses obtained via linear elasticity are more accurate than those computed by means of conventional plasticity theories \cite{IJSS2015,IJP2016,CM2017}. We assume small strains and define the strain tensor as,
\begin{equation}
    \bm{\epsilon}=\dfrac{1}{2}\left[\nabla \mathbf{u}^{T}+\nabla\mathbf{u}\right]
    \centering
\end{equation}

Accordingly, the Cauchy stress tensor can be derived from (\ref{Eq:FreeEnergy}) as,
\begin{equation}\label{Eq:StressConst}
\bm{\sigma}=\frac{\partial \Psi}{\partial \bm{\epsilon}} = g \left( \phi \right)  \left( \bm{C}_0 : \bm{\epsilon} \right) - K \bar{V}_H \left( C - C^0 \right) \bm{I}
\end{equation}
with the second term having a negligible effect in hydrogen embrittlement phenomena \cite{Hirth1980} and being subsequently disregarded in the analysis. As shown in (\ref{Eq:StressConst}), and discussed in Section \ref{Sec:PhaseFieldFracture}, the phase field model characterizes the loss in stiffness associated with bond breakage by means of the degradation function $g \left( \phi \right)$.

\subsubsection{Phase field fracture}

The scalar microstress $\omega$, work conjugate to the phase field $\phi$,
\begin{equation}
\omega=\frac{\partial \Psi}{\partial \phi} = - 2 (1-\phi) \, \psi_0  \left( \bm{\epsilon} \right) + G_c \left( \theta \right) \frac{1}{\ell} \phi
\end{equation}
and the vector microstress $\bm{\zeta}$, work conjugate to the phase field gradient $\nabla \phi$,
\begin{equation}
\bm{\zeta}=\frac{\partial \Psi}{\partial \nabla \phi} = G_c \left( \theta \right) \ell \nabla \phi 
\end{equation}
\noindent can be readily derived from the free energy (\ref{Eq:FreeEnergy}). Inserting these constitutive relations in the local balance (\ref{strong}b) renders (\ref{Eq:strong}) if the concentration gradient along the small region where $\nabla \phi \neq 0$ is neglected. On the one hand, the deformation of the solid drives the fracture phase field through the strain energy density $\psi_0$. On the other hand, the coupling between the phase field fracture and hydrogen transport is introduced by defining a hydrogen concentration-dependent crack surface energy that models decohesion enhancement by hydrogen. We provide this connection in a multi-scale fashion by presenting a quantum-mechanically informed phase field degradation law (see Section \ref{Sec:Henergydegradation}).

\subsubsection{Transport of diluted species}

The gradient of the chemical potential is the driving force for hydrogen diffusion. One can readily derive the chemical potential from the free energy (\ref{Eq:FreeEnergy}) by considering the relation between the occupancy and the number of sites, $\theta_L=C/N$,
\begin{equation}\label{Eq:ChePotential}
\mu=\frac{\partial \Psi}{\partial C} =\mu^0 + RT \ln \frac{\theta_L}{1 - \theta_L} - \bar{V}_H \sigma_H + G_c' \left( \theta \right) \theta' \left(\frac{1}{2 \ell} \phi^2 + \frac{\ell}{2}|\nabla \phi|^2 \right)
\end{equation}
and the constitutive framework is completed by the definition of the mass flux, which follows a linear Onsager relationship,
\begin{equation}\label{Eq:Flux}
    \bm{J}= - \frac{D C}{R T} \nabla \mu
\end{equation}
\noindent where $D$ is the diffusion coefficient. Thus, hydrogen atoms migrate from regions of high chemical potential to regions of low chemical potential. The coupling with the deformation of the solid is given by the stress-dependent part of $\mu$. As evident from (\ref{Eq:ChePotential}), hydrostatic tensile stresses increase hydrogen solubility in the lattice (lattice dilatation) by lowering the chemical potential. Correspondingly, the last term in (\ref{Eq:ChePotential}) enhances diffusion from cracked regions to pristine regions. However, the appropriate chemical boundary conditions in a propagating crack require careful consideration. First, one should note that the discussion is only relevant to experiments where the specimen is charged in-situ, as opposed to tests where hydrogen has been dissolved within the material prior to loading. Secondly, one would typically expect that the environment will promptly occupy the space created with crack advance. Consequently, a constant hydrogen concentration $C_{env}$ or, more appropriately, a constant chemical potential $\mu_{env}$, should be prescribed on a moving boundary $\partial B_{env}$ as dictated by the phase field $\phi$. In the finite element model this could be enforced by means of Dirichlet, Neumann or mixed boundary conditions \cite{Zhao2016,EFM2017}. However, we note in passing that crack growth rates $da/dt$ in hydrogen containing environments become load independent due to chemical reaction or mass transport limitation. A constant $\mu_{env}$ prescribed on a moving crack surface will fail to predict $da/dt$ measurements if a critical distance is not introduced in the failure criterion; such scheme could require a highly complex implementation. A physically sound framework could be key to quantitatively capture crack growth rates during the stage II plateau; since this regime has not been addressed in the present work, we assume for simplicity a negligible influence of the phase field on hydrogen diffusion.

\subsection{Numerical implementation}
\label{Sec:ABAQUS}

The finite element (FE) method is used to solve the coupled mechanical-diffusion-phase field problem. Using Voigt notation, the nodal values of the displacements, phase field and hydrogen concentration are interpolated as follows,
\begin{equation}\label{Eq:Discretization}
    \bm{u}=\sum_{i=1}^m \bm{N}_i \bm{u}_i \, , \,\,\,\,\,\,\,\,\,  \phi=\sum_{i=1}^m N_i \phi_i \, , \,\,\,\,\,\,\,\,\,  C=\sum_{i=1}^m N_i C_i
\end{equation}

\noindent where $m$ is the number of nodes and $\bm{N}_i$ are the interpolation matrices - diagonal matrices with the nodal shape functions $N_i$ as components. Accordingly, the corresponding gradient quantities can be discretized by,
\begin{equation}
    \bm{\varepsilon}=\sum_{i=1}^m \bm{B}^u_i \bm{u}_i \, , \,\,\,\,\,\,\,\,\,  \nabla \phi=\sum_{i=1}^m \bm{B}_i \phi_i \, , \,\,\,\,\,\,\,\,\,  \nabla C=\sum_{i=1}^m \bm{B}_i C_i
\end{equation}

\noindent Here, $\bm{B}_i$ are vectors with the spatial derivatives of the shape functions and $\bm{B}^u_i$ denotes the standard strain-displacement matrices.

\subsubsection{FE discretization of the deformation-phase field problem}
The weak form for the deformation problem can be written, in its most general form, as
\begin{equation}\label{Eq:WeakDisp}
\int_{\Omega} \left( g \left( \phi \right) \bm{\sigma}_0 : \delta \bm{\epsilon} - \mathbf{b} \cdot \delta \mathbf{u} \right) \, \mathrm{d}V - \int_{\partial \Omega_{h}} \left( \mathbf{h} \cdot \delta \mathbf{u} \right) \, \mathrm{d}S = 0
\end{equation}
where $\bm{\sigma}_0$ is the Cauchy stress tensor of the undamaged solid. Making use of the finite element discretization outlined above and considering that (\ref{Eq:WeakDisp}) must hold for arbitrary values of $\delta \bm{u}$, the discrete equation corresponding to the equilibrium condition can be expressed as the following residual with respect to the displacement field,
\begin{equation} \label{residualStagPhi}
    \bm{r}_{i}^u =\int_\Omega \left[(1-\phi)^{2}+k\right] {(\bm{B}_{i}^u)}^{T} \bm{\sigma_{0}} \, \mathrm{d}V - \int_\Omega \bm{N}_i^T \bm{b} \, \mathrm{d}V - \int_{\partial \Omega_h} \bm{N}_i^T \bm{h} \, \mathrm{d}S   
\end{equation}

Similarly, the out-of-balance force residual with respect to the evolution of the crack phase field is obtained by discretizing (\ref{Eq:weak}) and considering the contribution from the scalar phase field microtraction (\ref{Eq:PVW}),
\begin{equation}
    r_{i}^{\phi}= \int_\Omega \left[ -2(1-\phi) N_{i} \, H + G_c \left( \theta \right) \left(\dfrac{1}{\ell} N_{i} \, \phi
    + \ell \bm{B}_{i}^T \nabla \phi \right) \right] \, \mathrm{d}V - \int_{\partial \Omega_f} N_{i} f \, \text{d}S
\end{equation}
\noindent where we have introduced the so-called history variable field \textit{H} to ensure irreversibility,
\begin{equation} \label{history}
    H =
    \begin{cases}
      \psi_{0}(\bm{\epsilon})  & \quad \text{if} \hspace{2mm} \psi_{0}(\bm{\epsilon}) > H_{t} \\
      H_{t}  & \quad \text{otherwise}\\
    \end{cases}
    \centering
\end{equation}
\noindent Here, \(H_{t}\) is the previously calculated energy at time increment \(t\). Thus, the history field satisfies the Kuhn-Tucker conditions. The components of the consistent stiffness matrices can be obtained by differentiating the residuals with respect to the incremental nodal variables:
\begin{equation}\label{Eq:Ku}
    \bm{K}_{ij}^{\bm{u}} = \frac{\partial \bm{r}_{i}^{\bm{u}}}{\partial \bm{u}_{j}} = \int_\Omega \left[(1-\phi)^2+ k\right] {(\bm{B}_i^{\bm{u}})}^T \bm{C}_0 \bm{B}_j^{\bm{u}} \, \mathrm{d}V  
\end{equation}
\begin{equation}
    \bm{K}_{ij}^\phi = \dfrac{\partial r_{i}^{\phi}}{\partial \phi_{j}} = \int_\Omega \left[ \left( 2H + \dfrac{G_c \left( \theta \right)}{\ell} \right) N_{i} N_{j} + G_c \left( \theta \right) \ell \bm{B}_i^T \bm{B}_j \right] \, \mathrm{d}V    
\end{equation}

\subsubsection{FE discretization of mass transport}

Substituting the constitutive equation for the chemical potential (\ref{Eq:ChePotential}) into the Fick-type relation for the mass flux (\ref{Eq:Flux}), one reaches
\begin{equation}
   \bm{J}=- \frac{DC}{\left(1 - \theta_L \right)} \left(\frac{\nabla C}{C} - \frac{\nabla N}{N} \right) + \frac{D}{RT} C \bar{V}_H \nabla \sigma_H
\end{equation}
\noindent which, after making the common assumptions of low occupancy ($\theta_L \ll 1$) and constant interstitial sites concentration ($\nabla N=0$), renders
\begin{equation}
   \bm{J}=-D \nabla C + \frac{D}{RT} C \bar{V}_H \nabla \sigma_H
\end{equation}
Substituting in (\ref{Eq:PVW}), the hydrogen transport equation becomes
\begin{equation}\label{Eq:WeakC}
    \int_\Omega \left[ \delta C \left( \frac{1}{D} \frac{dC}{dt} \right) + \nabla \delta C  \nabla C  - \nabla \delta C \left( \frac{\bar{V}_H C}{RT}  \nabla \sigma_H \right) \right] \, \mathrm{d}V = - \frac{1}{D} \int_{\partial \Omega_q} \delta C q \, \mathrm{d}S
\end{equation}

A residual vector can be readily obtained by discretizing (\ref{Eq:WeakC}), given that $\delta C$ indicates an arbitrary virtual variation of the hydrogen concentration: 
\begin{equation}
     r_{i}^C= \int_\Omega \left[ N_i^T \left( \frac{1}{D} \frac{dC}{dt} \right) + \bm{B}_i^T  \nabla C  - \bm{B}_i^T \left( \frac{\bar{V}_H C}{RT}  \nabla \sigma_H \right)  \right] \, \mathrm{d}V + \frac{1}{D} \int_{\partial \Omega_q} N_i^T q \, \mathrm{d}S
\end{equation}

From which a diffusivity matrix can be defined,
\begin{equation}
    \bm{K}_{ij}^C = \int_\Omega \left( \bm{B}_i^T  \bm{B}_j  - \bm{B}_i^T \frac{\bar{V}_H}{RT}   \nabla \sigma_H N_j   \right) \, \mathrm{d}V    
\end{equation}

\noindent where the discretization given in Eq. (\ref{Eq:Discretization}) has also been employed to interpolate the time derivatives of the nodal concentrations. The diffusivity matrix is affected by the gradient of the hydrostatic stress, $\sigma_H$, which is computed at the integration points from the nodal displacements, extrapolated to the nodes by means of the shape functions, and subsequently multiplied by $\bm{B}$ to compute $\nabla \sigma_H$. We employ 8-node quadrilateral finite elements with $C_0$ continuity.\\

At the same time, one can readily identify a concentration capacity matrix,
\begin{equation}
    \bm{M}_{ij} = \int_\Omega N_i^T \frac{1}{D} N_j \, \mathrm{d}V        
\end{equation}

\noindent and a diffusion flux vector,
\begin{equation}
   \bm{F}_{i} = - \frac{1}{D} \int_{\partial \Omega_q} N_i^T q \,\, \textnormal{d}S
\end{equation}

\noindent Accordingly, the discretized hydrogen transport equation reads,
\begin{equation}
    \bm{K}^C \bm{C} + \bm{M} \dot{\bm{C}} = \bm{F}
\end{equation}

\subsubsection{Coupled scheme}

The deformation, diffusion and phase field fracture problems are weakly coupled. First, mechanical deformation impacts diffusion through the stress field, governing the pressure dependence of the bulk chemical potential. Secondly, mass transport affects the fracture resistance via hydrogen buildup in the fracture process zone, reducing the critical energy release rate. And thirdly, the hydrogen-sensitive phase field degrades the strain energy density of the solid.\\

We solve the linearized finite element system,
\begin{equation}
\begin{bmatrix}
  \bm{K}^u & 0 & 0\\
  0 & \bm{K}^\phi & 0 \\
  0 & 0 & \bm{K}^C 
 \end{bmatrix} \begin{bmatrix} \bm{u} \\ \bm{\phi} \\ \bm{C} \end{bmatrix} + 
 \begin{bmatrix}
  0 & 0 & 0\\
  0 & 0 & 0 \\
  0 & 0 & \bm{M}
 \end{bmatrix} \begin{bmatrix} \dot{\bm{u}} \\ \dot{\bm{\phi}} \\ \dot{\bm{C}} \end{bmatrix}= \begin{bmatrix} \bm{r}^u \\ \bm{r}^\phi \\ \bm{r}^C \end{bmatrix}
\end{equation}

\noindent by means of a time parametrization and an incremental-iterative scheme in conjunction with the Newton-Raphson method. Due to its robustness, the staggered solution scheme proposed by Miehe et al. \cite{Miehe2010a} is adopted to solve the deformation-phase field coupling; a time increment sensitivity analysis is conducted in all computations.\\

The modeling framework is implemented in the commercial finite element package ABAQUS via a user element subroutine. The code can be downloaded from www.empaneda.com/codes, including documentation describing the details of the implementation in ABAQUS, verification case studies and integration strategies. Post-processing of the results is carried out by means of Abaqus2Matlab \cite{AES2017}.

\section{Results}
\label{Sec:Results}

A number of case studies of particular interest are addressed to show the capabilities of the modeling framework presented. Firstly, we model the benchmark problem of a cracked square plate subjected to tension, verifying our implementation with the results by Miehe et al. \cite{Miehe2010a}, and we subsequently assess the role of hydrogen (Section \ref{Sec:MiehePlate}). Then, in Section \ref{Sec:AxisymmetricBars}, hydrogen assisted failure in notched components is investigated, and model predictions are compared with the experimental results of Wang et al. \cite{Wang2005}. The experimental work by Olden et al. \cite{Olden2009} is reproduced in Section \ref{Sec:Olden}, where the model is used to estimate cracking thresholds under constant loading conditions. We investigate internal hydrogen assisted cracking and compare with the compact tension experiments by Thomas et al. \cite{Thomas2003} in Section \ref{Sec:CompactTension}. Finally, Section \ref{Sec:PittingCorrosion} deals with complex crack propagation from defects that typically arise due to pitting corrosion.

\subsection{Cracked square plate subjected to tension in a hydrogenous environment}
\label{Sec:MiehePlate}

\begin{figure}[H] 
    \centering
    \includegraphics[scale=1.5]{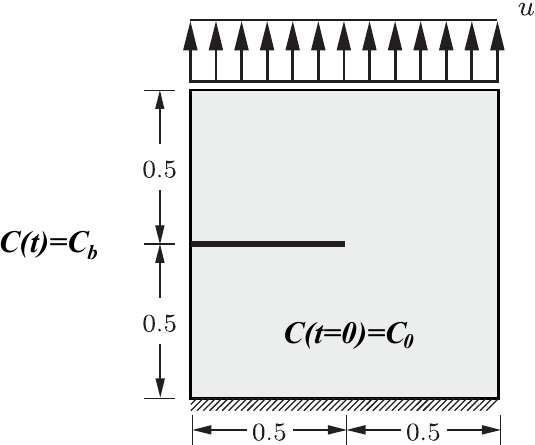}
    \caption{Schematic description of the geometry and configuration of the cracked square plate case study.}
    \label{fig:PlateBC}
\end{figure}

We model a square plate with a horizontal crack going from the left side to the center of the specimen. The geometric setup, dimensions and boundary conditions are given in Fig. \ref{fig:PlateBC}. We load the plate by prescribing the vertical displacement in the upper edge, and restrict both vertical and horizontal displacements in the lower boundary. A uniform distribution of hydrogen concentration is assigned throughout the specimen as initial boundary condition: $C(t=0)=C_0$. In addition, a constant hydrogen concentration $C=C_b$ is prescribed at the boundaries during the entire numerical experiment. All the outer boundaries of the specimen, including the crack faces, are in contact with the environment. These boundary conditions mimic the usual scenario in field applications and laboratory experiments, where the specimen is pre-charged before loading in the same aggressive environment; accordingly, we assume $C_b=C_0$. For materials with high hydrogen diffusivity, the use of generalized Neumann-type boundary conditions \cite{Turnbull2015} or $\sigma_H$-dependent Dirichlet boundary conditions \cite{DiLeo2013,IJHE2016,Diaz2016b} could provide a more accurate description of the electrochemistry-diffusion interface. As in \cite{Miehe2010a}, we adopt the following material properties: Young's modulus $E=210$ GPa, Poisson's ratio $\nu=0.3$, and critical energy release rate $G_{c} \left( 0 \right)=2.7$ MPa mm. We assume an iron-based material and consequently adopt a hydrogen damage coefficient of $\chi=0.89$ and a partial molar volume of $\bar{V}_H=2000$ mm$^3$/mol \cite{Hirth1980}. For this benchmark study, we assume that the load is prescribed at a sufficiently low rate such that hydrogen equilibrium can be achieved. Thus, from the diffusion coefficient and the test time $t_f$ one can define a diffusion length $\ell_C = \sqrt{D t_f}$, which is taken to be significantly larger than the characteristic length of the fracture problem ($\ell_C >> \ell$). A very fine mesh is employed in the expected crack propagation area after a sensitivity study, ensuring that the characteristic element length $h$ is sufficiently small to resolve the fracture process zone ($h< \ell/7.5$). Namely, the square plate is discretized by means of 43110 plane strain quadratic elements.

\begin{figure}[H]
\centering
\includegraphics[scale=1.1]{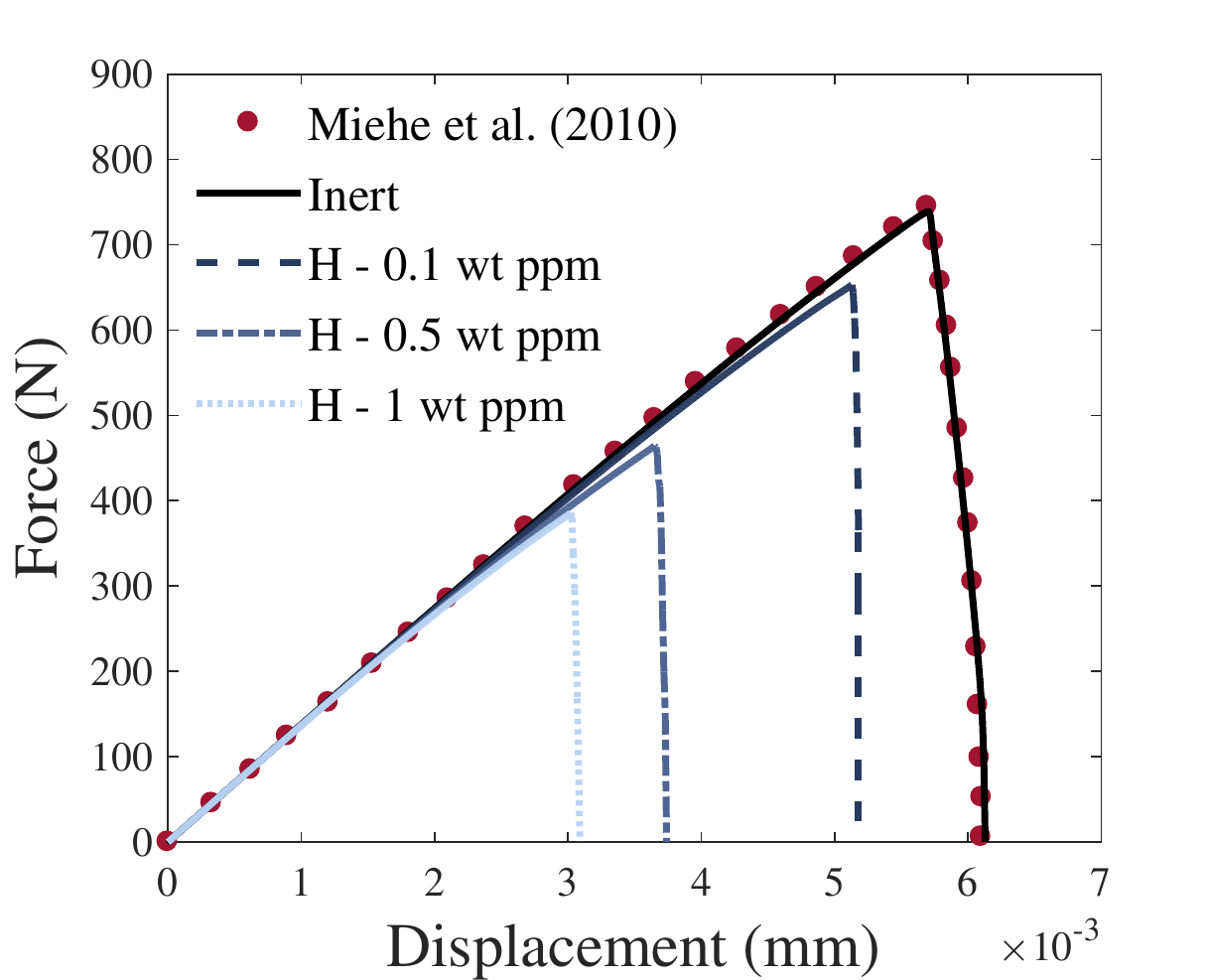}
\caption{Cracked square plate subjected to tension. Load-deflection curves for different environments. The computations for $C_0=C_b=0$ wt ppm match the results by Miehe et al. \cite{Miehe2010a} (symbols). Material properties: $E=210$ GPa, $\nu=0.3$, $G_c=2.7$ KJ/m$^2$, $\ell=0.0075$ mm, $\bar{V}_H=2000$ mm$^3$/mol, $\ell_C >> \ell$, and $T=300$ K.}
\label{fig:LoadDeflectionMiehe}
\end{figure}

To assess the capabilities of the model in capturing hydrogen-induced damage we compute the mechanical response under three different environments that range from 1 wt ppm, which corresponds to a 3\% NaCl aqueous solution, to 0.1 wt ppm. Fig. \ref{fig:LoadDeflectionMiehe} shows the results computed in terms of the load-displacement curve. As shown in the figure, the model shows a strong sensitivity to the environment, with the peak load decreasing with the hydrogen content, in agreement with experimental observations. The numerical results in the absence of hydrogen quantitatively agree with those by Miehe et al. \cite{Miehe2010a}, verifying the numerical implementation.

\begin{figure}[H] 
    \centering
    \includegraphics[scale=0.4]{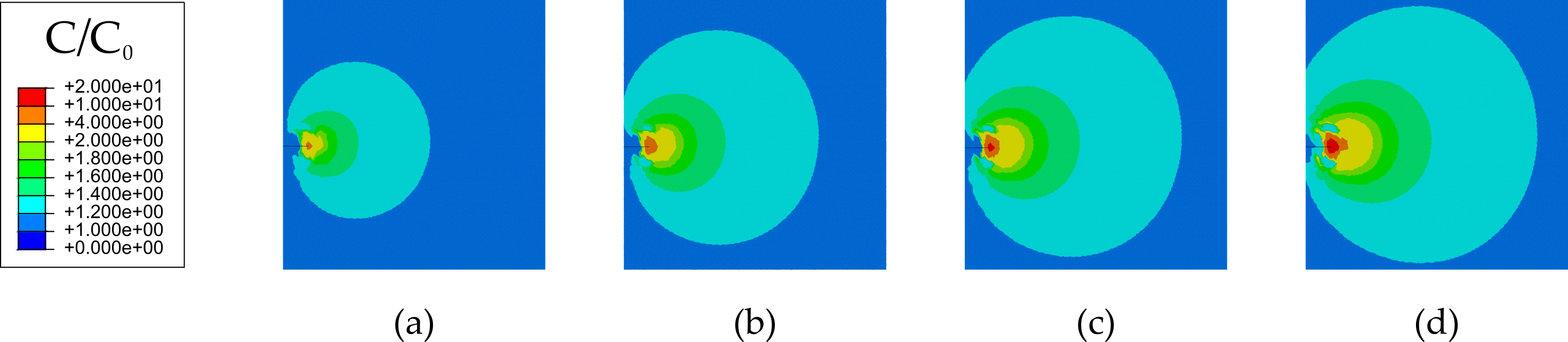}
    \caption{Hydrogen buildup in the fracture process zone. Contours at (a) $u=0.0015$ mm, (b) $u=0.0019$ mm, (c) $u=0.0021$ mm, and (d) $u=0.0023$ mm for $C_0=0.5$ wt ppm.}
    \label{fig:ContoursH}
\end{figure}

As shown in Fig. \ref{fig:ContoursH}, hydrogen accumulates in the vicinity of the crack as the external load and the hydrostatic stress increase.

\begin{figure}[H] 
    \centering
    \includegraphics[scale=0.45]{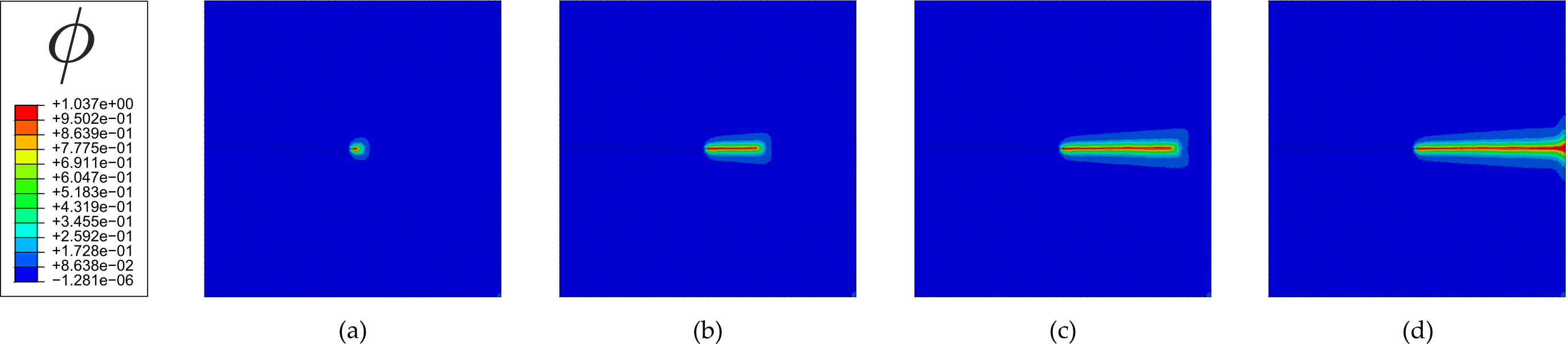}
    \caption{Crack propagation contours at (a) $u=0.0057$ mm, (b) $u=0.0059$ mm, (c) $u=0.0061$ mm, and (d) $u=0.0062$ mm for $C_0=0$ wt ppm.}
    \label{fig:ContoursPhi}
\end{figure}

Atomic hydrogen concentration lowers the critical energy release rate of the material, facilitating crack initiation and subsequent propagation. Fig. \ref{fig:ContoursPhi} shows the resulting fracture patterns at different load steps with the red color denoting the case where $\phi=1$. The trend depicted by the crack path is the same in all the environments considered. 

\subsection{Hydrogen assisted failure in notched cylindrical bars}
\label{Sec:AxisymmetricBars}

\begin{figure}[H] 
    \centering
    \includegraphics[scale=0.9]{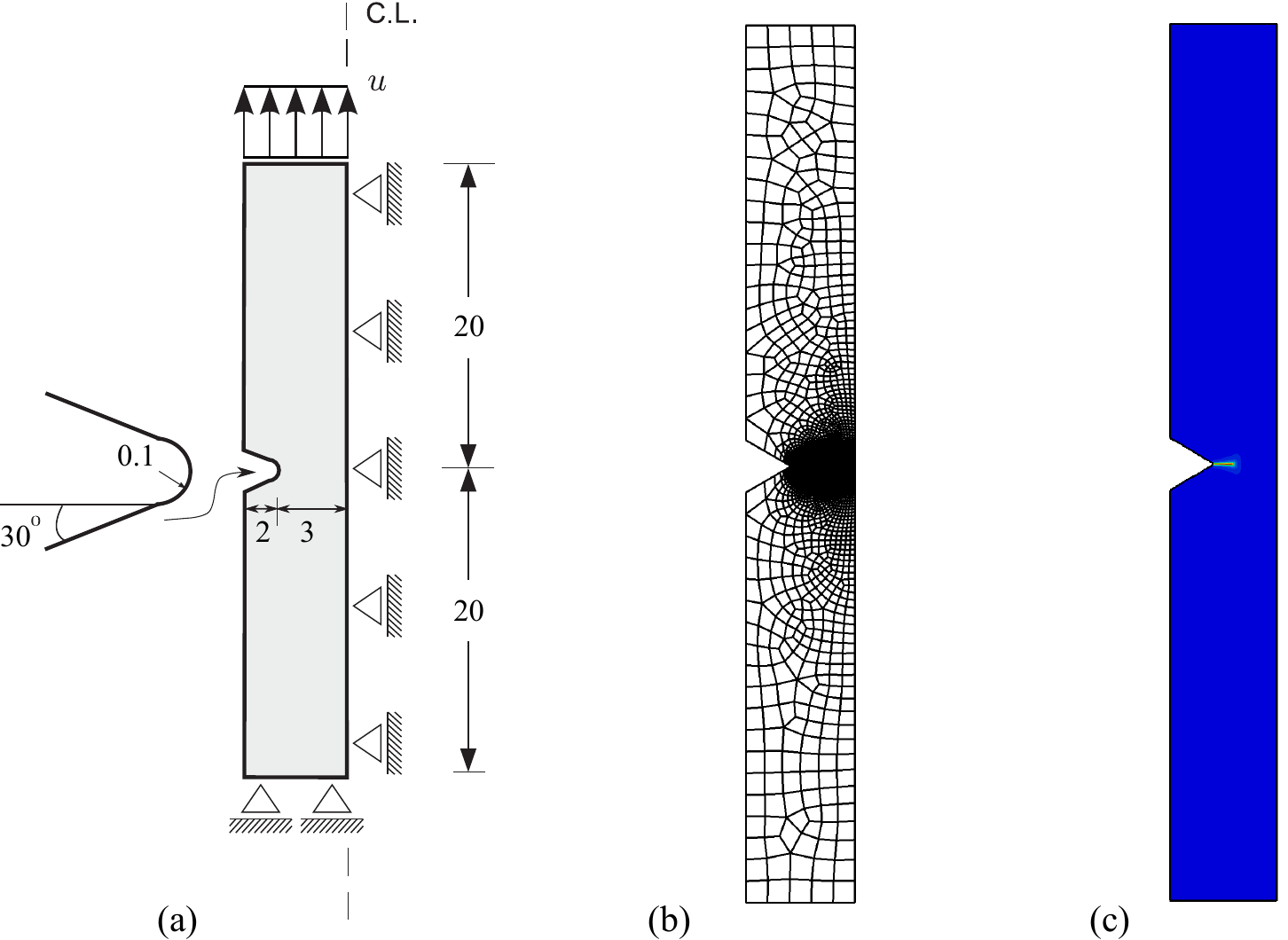}
    \caption{Schematic description of the notch cylindrical bar: (a) geometry, (b) mesh and (c) crack path when the remote stress equals $\sigma_f$. Dimensions are given in mm.}
    \label{fig:Bar}
\end{figure}

We now turn our attention to hydrogen assisted failure from notched specimens under rising load and different environmental conditions. We reproduce the experimental work by Wang et al. \cite{Wang2005} on cylindrical bars. The plane strain implementation can be readily modified to deal with axisymmetric samples by changing the strain-displacement matrix and integrating the discretized equations in polar coordinates. Fig. \ref{fig:Bar} shows the specimens geometry, the mesh employed and the path depicted by the crack that initiates at the notch tip. A total of 52955 axisymmetric 8-node elements are used to discretize the bar. The length of the characteristic element is chosen to be 10 times smaller than $\ell$ to assure mesh-independent results and resolve the fracture process zone.\\

\begin{figure}[H] 
    \centering
    \includegraphics[scale=1]{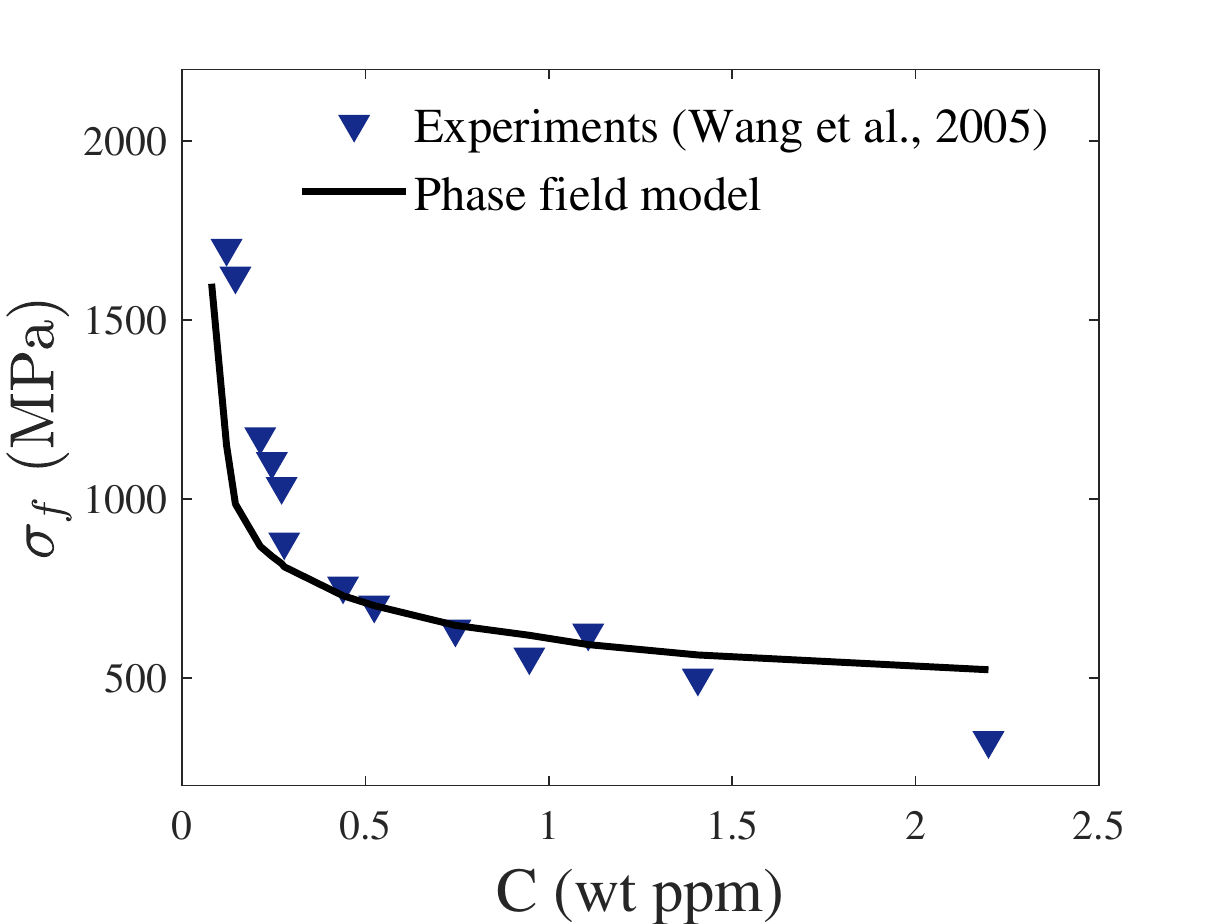}
    \caption{Notched bar subjected to tension test. Net section strength $\sigma_f$ as a function of the hydrogen concentration $C$ in AISI 4135 steel. Material properties: $E=210$ GPa, $\nu=0.3$, $\bar{V}_H=2000$ mm$^3$/mol, $D=3.8 \times 10^{-5}$ mm$^2$/s, $\ell=0.015$ mm, and $G_c=64$ KJ/m$^2$.}
    \label{fig:Wang}
\end{figure}

We mimic the conditions of the experiments on AISI 4135 steel bars by Wang et al. \cite{Wang2005}. Uniaxial tension tests are conducted at a nominal strain rate of $\dot{\varepsilon}=8.3 \times 10^{-7}$ s$^{-1}$. The experiments are conducted at room temperature. The material properties for AISI 4135 steel are: Young's modulus $E=210$ GPa, Poisson's ratio $\nu=0.3$, partial molar volume of $\bar{V}_H=2000$ mm$^3$/mol, and diffusion coefficient $D=3.8 \times 10^{-5}$ mm$^2$/s. The specimens are precharged and subsequently loaded in the same environment; hydrogen concentrations range from 0.08 to 2.2 wt ppm. For each case, the maximum net section stress is computed at the remote boundary. Fig. \ref{fig:Wang} shows the failure stress obtained in each numerical experiment as a function of the hydrogen content. The results from Wang et al. \cite{Wang2005} are also shown; a critical energy release rate of $G_c=64$ KJ/m$^2$ provided the best fit to the data.

Phase field predictions reveal a good agreement with the experimental data. Increasing the hydrogen content translates into a reduction of the critical energy release rate, which in turn lowers the strength. Hydrogen accumulates in the vicinity of the notch tip, where tensile hydrostatic stresses are larger. The strength changes drastically within low hydrogen levels but eventually saturates: increasing the hydrogen content beyond 0.5 wt ppm leads to minor changes in the failure stress. The saturation stage is captured due to the presence of the Langmuir-McLean isotherm (\ref{Eq:Langmuir-McLean}). The embrittlement law seems to capture the experimentally observed response well. Discrepancies are larger in the experiments with very low hydrogen concentrations, where grain boundary decohesion may not be the only mechanism assisting failure.

\subsection{Cracking threshold for duplex steels in sea water under constant load}
\label{Sec:Olden}

\begin{figure}[H] 
  \makebox[\textwidth][c]{\includegraphics[width=1.2\textwidth]{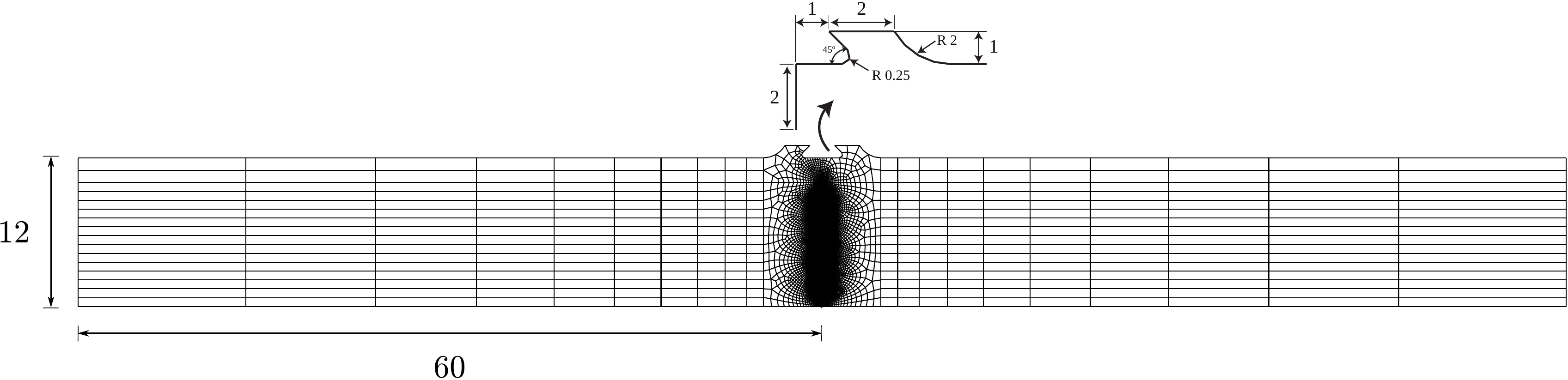}}%
  \caption{Geometry and mesh of the single edge notched specimen. Dimensions are given in mm.}
  \label{fig:SENTgeometry}
\end{figure}

We investigate the capabilities of the present model to capture hydrogen damage under constant mechanical load by reproducing the experiments by Olden and co-workers \cite{Olden2009}. In their work, single-edge notched tensile (SENT) specimens were immersed in slowly circulating artificial sea water at 4$^\circ$ C under cathodic protection conditions. The material under consideration is a 25\%Cr duplex stainless steel with $E=200$ GPa and $\nu=0.3$; despite its wide use in sub-sea applications, duplex steels are sensitive to hydrogen embrittlement at low corrosion protection potentials \cite{Olden2009}. The geometry of the fatigue pre-cracked SENT specimens is shown in Fig. \ref{fig:SENTgeometry}, along with the finite element mesh employed. A total of 98748 plane strain quadratic elements are used in the discretization. A constant net stress is prescribed at the remote boundaries and a subsurface concentration of 1 wt ppm is assumed in the outer surfaces, based on the experimental conditions (3.5\% NaCl solution and a corrosion protection potential of $-1050$ mV$_{SCE}$). The tests have been conducted at room temperature and the diffusion coefficient has been measured as $D=3.7 \times 10^{-6}$ mm$^2$/s.

The initial mechanical load introduces a certain amount of damage that is insufficient to propagate the crack. The magnitude of the phase field variable increases and eventually triggers crack growth as time goes by and hydrogen diffuses towards the fracture process zone. Since the mechanical load is kept constant, damage is caused by the reduction in the critical energy release rate with time due to hydrogen diffusion.\\

\begin{figure}[H] 
    \centering
    \includegraphics[scale=1]{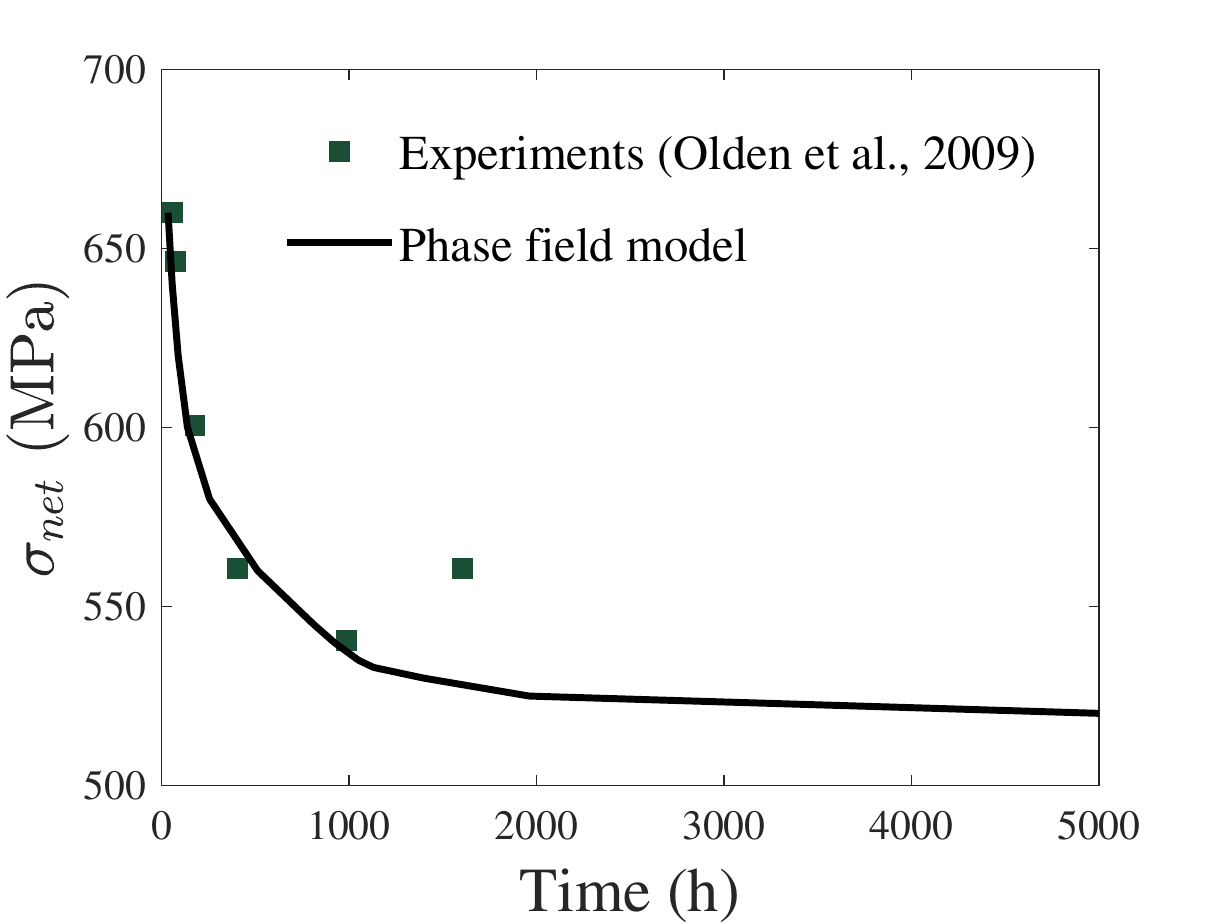}
    \caption{SENT specimen subjected to constant load. Net section stress $\sigma_{net}$ as a function of time in duplex stainless steel. Material properties: $E=200$ GPa, $\nu=0.3$, $\bar{V}_H=2000$ mm$^3$/mol, $D=3.7 \times 10^{-6}$ mm$^2$/s, $\ell=0.06$ mm, and $G_c=27.4$ KJ/m$^2$.}
    \label{fig:Olden}
\end{figure}

We compare our cracking threshold predictions with the experimental data by Olden et al. \cite{Olden2009}. Results are shown in Fig. \ref{fig:Olden} in terms of net stress versus time. As shown in the figure, we are able to reach a good agreement with the tests and extend our predictions over time scales that are impractical in laboratory experiments. This allows us to identify the net stress level below which cracking will not occur: $\sigma_{TH}=510$ MPa. Thus, safe regimes of operation can be efficiently identified by means of the present phase field formulation. However, one should note that the lack of experimental data at larger time scales introduces a degree of uncertainty in the numerical predictions of this case study. Specifically, more tests within the range of 1000-2000 h would have been desirable given the scatter observed for net stresses of 540-560 MPa.\\

\subsection{Internal hydrogen assisted cracking in compact tension specimens}
\label{Sec:CompactTension}

Hydrogen assisted cracking frequently takes place from hydrogen dissolved within the material prior to loading. This phenomenon is termed internal hydrogen assisted cracking, and atomic hydrogen is typically introduced by manufacturing operations, such as casting, welding, electrochemical machining, or heat treatments. Mechanical loading redistributes the internal hydrogen content to the fracture process zone and promotes cracking.\\ 

\begin{figure}[H]
\makebox[\linewidth][c]{%
        \begin{subfigure}[h]{0.4\textwidth}
                \centering
                \includegraphics[scale=0.8]{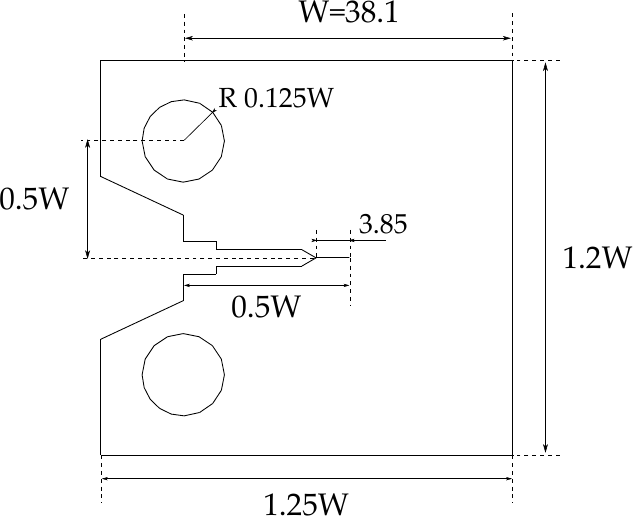}
                \caption{}
                \label{fig:Geometry}
        \end{subfigure}
        \begin{subfigure}[h]{0.4\textwidth}
                \centering
                \includegraphics[scale=0.44]{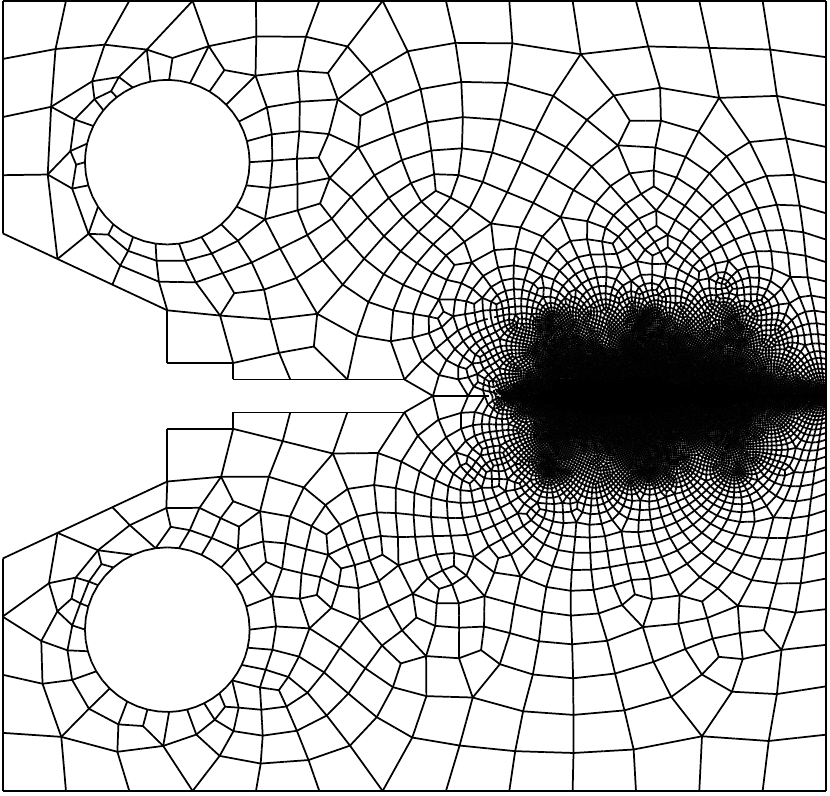}
                \caption{}
                \label{fig:Mesh}
        \end{subfigure}
        \begin{subfigure}[h]{0.4\textwidth}
                \centering
                \includegraphics[scale=0.07]{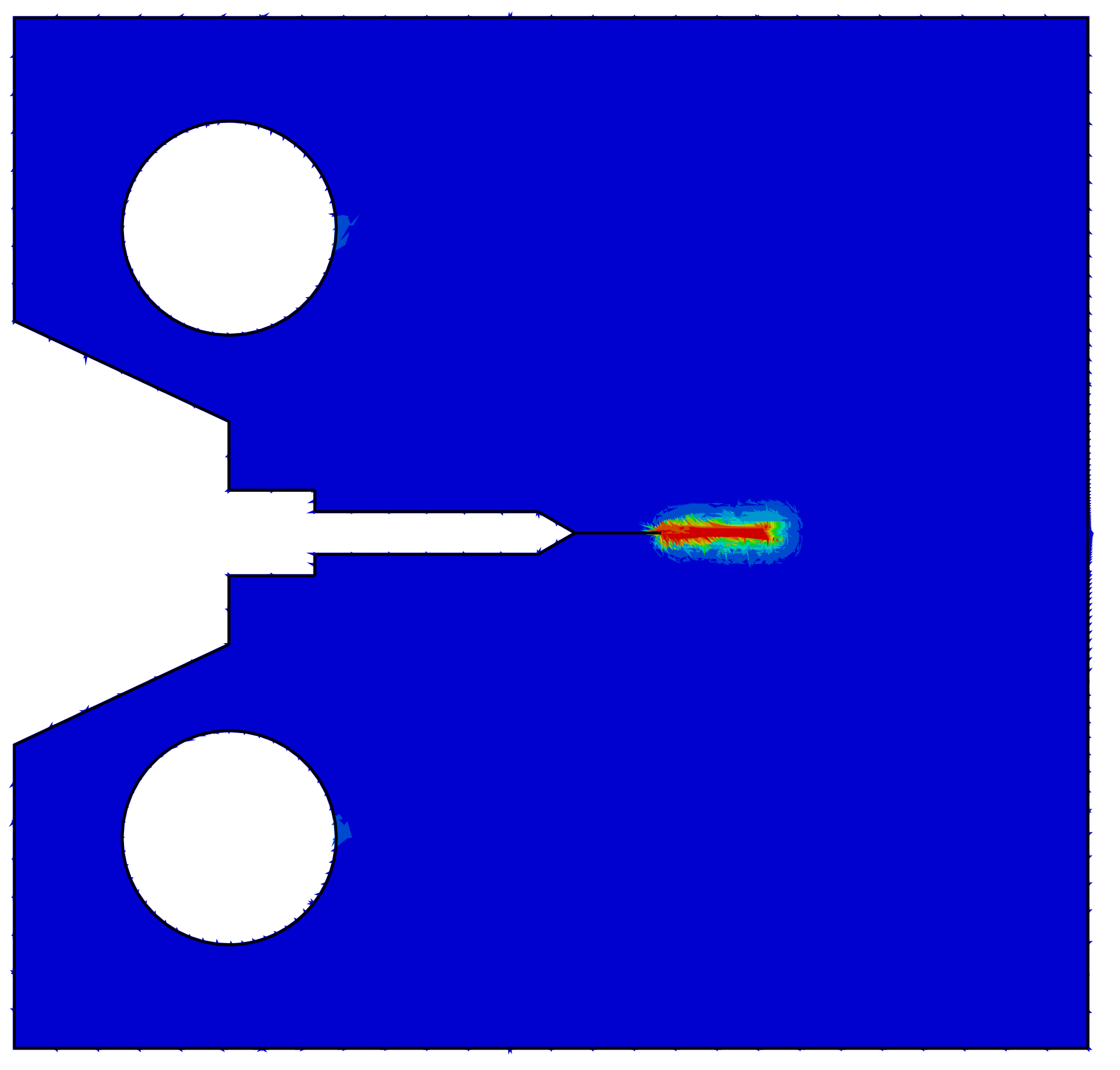}
                \caption{}
                \label{fig:Crack}
        \end{subfigure}}        
       
        \caption{Schematic description of the Compact Tension experiment: (a) geometry, (b) mesh and (c) crack path. Dimensions are given in mm.}\label{fig:CT}
\end{figure}

We show the capabilities of the present phase field formulation to model internal hydrogen assisted cracking by mimicking the experiments by Thomas et al. \cite{Thomas2003} on AerMet100. AerMet100 is secondary-hardening quenched and tempered martensitic steel developed for high-performance aerospace applications. An electroplated coating is used for corrosion resistance, but atomic hydrogen can be codeposited in the plating and steel substrate \cite{Thomas2003}. In this case study, we model pre-cracked Compact Tension specimens that have been pre-charged prior to loading. Our aim is to predict the threshold stress intensity for internal hydrogen embrittlement $K_{TH}$ as a function of the hydrogen content. Fig. \ref{fig:Geometry} shows the geometry and dimensions of the specimens, following \cite{Thomas2003}, and in agreement with the ASTM E 1820 Standard. Fig. \ref{fig:CT} also shows the mesh (Fig. \ref{fig:Mesh}) and the characteristic mode I crack propagation path (Fig. \ref{fig:Crack}). A total of 34792 quadratic plane strain elements have been used and the stress intensity factor at crack initiation is computed following the ASTM E 1820 Standard. The material properties related to the mechanical response are Young's modulus $E=194400$ MPa, and Poisson's ratio $\nu=0.3$. The experiments have been conducted at room temperature and the diffusion coefficient for AerMet100 has been measured as $D=2 \times 10^{-4}$ mm$^2$/s \cite{Thomas2003}. Mimicking the experiments, a loading rate of $\dot{K}=2.2 \times 10^{-4}$ MPa$\sqrt{m}$/s is employed. The length parameter is chosen to be $\ell=0.3$ mm; 7 times larger than the characteristic element length along the extended crack plane.

\begin{figure}[H] 
    \centering
    \includegraphics[scale=1]{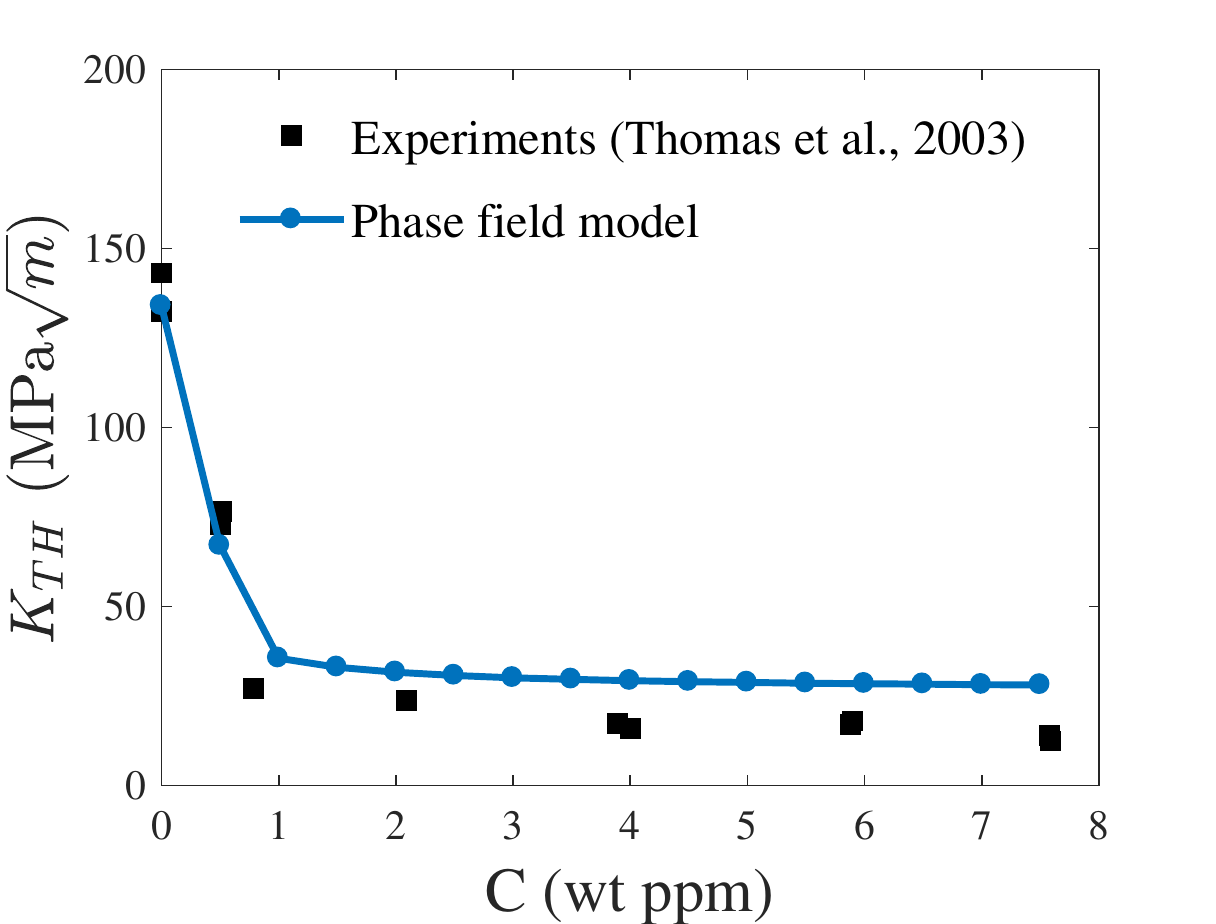}
    \caption{Internal hydrogen cracking in Compact Tension specimens. Cracking threshold $K_{TH}$ as a function of the hydrogen concentration $C$ in AerMet100. Material properties: $E=194400$ MPa, $\nu=0.3$, $\bar{V}_H=2000$ mm$^3$/mol, $D=2 \times 10^{-4}$ mm$^2$/s, $\ell=0.3$ mm, and $G_c=30$ KJ/m$^2$.}
    \label{fig:Thomas}
\end{figure}

Fig. \ref{fig:Thomas} shows the results obtained in terms of $K_{TH}$ versus the uniform hydrogen concentration attained prior to loading. A value of $G_c=30$ KJ/m$^2$ is assigned to the critical energy release rate to fit the experiments in air. As shown in Fig. \ref{fig:Thomas}, the numerical predictions qualitatively capture the experimental trends. Namely, the cracking threshold decreases significantly with increasing hydrogen content. Internal hydrogen embrittlement is particularly severe for this very high strength steel, with $K_{TH}$ being roughly a tenth of the fracture toughness in inert environments. While the agreement with experiments is considered satisfactory, finite element results fail to quantitatively predict the degree of embrittlement at large hydrogen concentrations. Plausibly, the differences are a consequence of fitting $G_c$ with stress intensity measurements that pertain to fracture processes where nucleation, growth and coalescence of microvoids are the mechanisms governing fracture. The value assigned to the trapping energy at grain boundaries may constitute another source of discrepancy; accurate measurements of trap characteristics remain a challenging task.

\subsection{Crack coalescence from corrosion pits}
\label{Sec:PittingCorrosion}

Hydrogen assisted cracking frequently manifests in combination with other corrosion damaging effects. For example, this is often the case in aqueous solutions where cathodic protection has been used to enhance corrosion resistance \cite{Sanchez2017}. Cracking in pipelines and other structural components often commences at the sites of corrosion pits, and is accelerated by hydrogen embrittlement. We investigate the capabilities of the present phase field formulation to predict complex crack propagation due to existing defects. This capability could be a game-changer in the energy industry, where structural integrity assessment often involves estimating crack growth and coalescence from well-mapped corrosion pits as a function of time.\\

\begin{figure}[H] 
    \centering
    \includegraphics[scale=0.5]{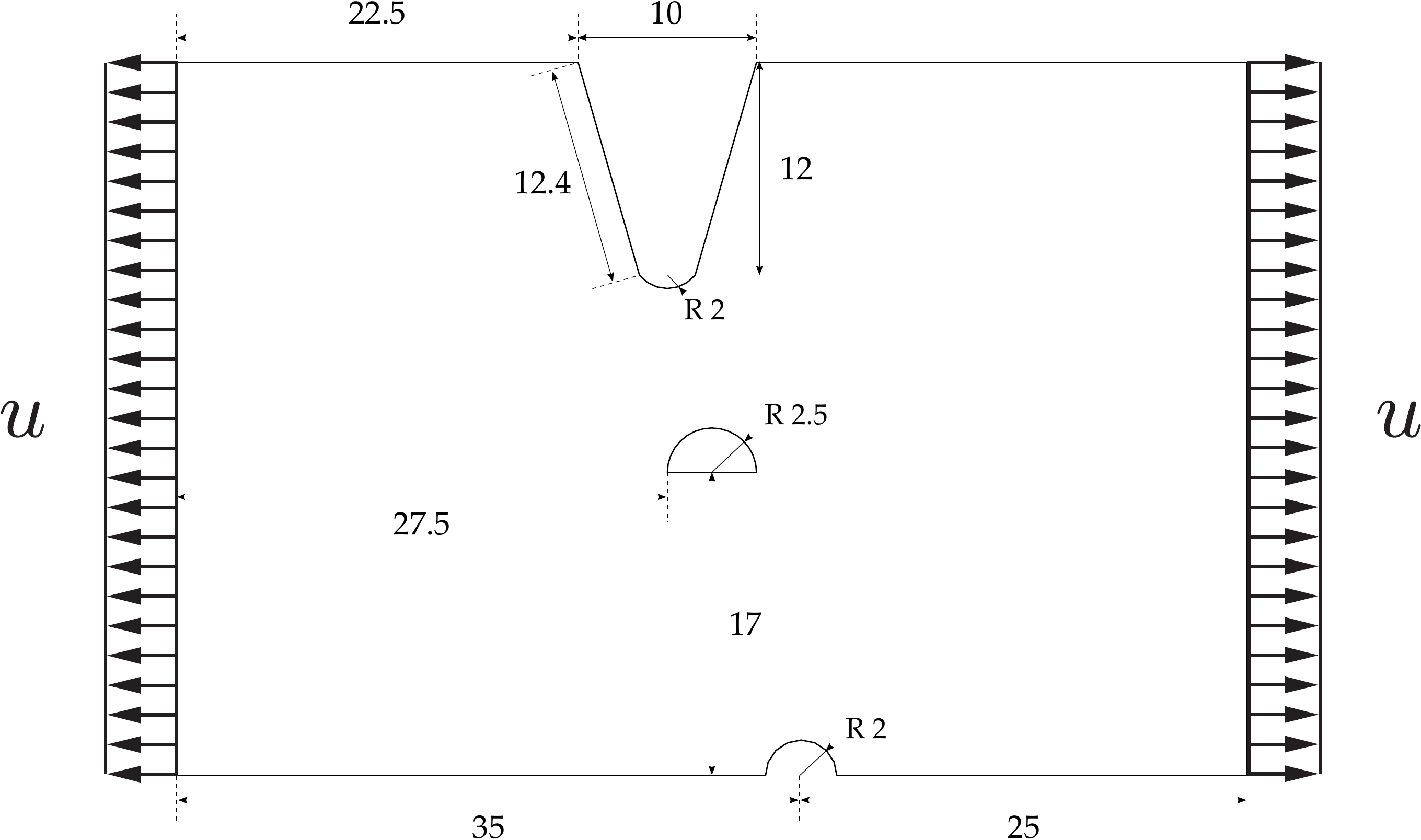}
    \caption{Geometry and loading configuration of a plate with existing defects, dimensions are given in mm. The defects are introduced in the initial rectangular geometry by prescribing the initial condition $\phi=1$.}
    \label{fig:GeometryPit}
\end{figure}

We consider a plate with three preexisting defects, as outlined in Fig. \ref{fig:GeometryPit}. The defects, which typically originate due to anodic dissolution, are defined in the numerical model by prescribing an initial history field $\phi=1$. Introducing pre-existing cracks as heterogeneities in the phase field order parameter (as opposed to discrete geometrical discontinuities) has proven to yield predictions closer to analytic fracture mechanics solutions \cite{Klinsmann2015}. We assume $E=200$ GPa, $\nu=0.3$, and energy release rate $G_c=90$ KJ/m$^2$. The length parameter is 6 times larger than the characteristic element size, $\ell=0.6$ mm. A mesh of 30908 quadratic quadrilateral plane strain elements is employed to model the plate. We consider a diffusion coefficient on the order of that of fcc iron at room temperature, $D=1 \times 10^{-8}$ mm$^2$/s, and prescribe an initial hydrogen concentration of $C_0=1$ wt ppm throughout the specimen. This hydrogen concentration is maintained during the test at the sides of the specimen and in the corrosion pits. As outlined in Fig. \ref{fig:GeometryPit}, we prescribe a uniform displacement at the sides, being the loading rate $\dot{u}=0.0416$ $\mu$m/s.

\begin{figure}[H]
\makebox[\linewidth][c]{%
        \begin{subfigure}[h]{0.6\textwidth}
                \centering
                \includegraphics[scale=0.4]{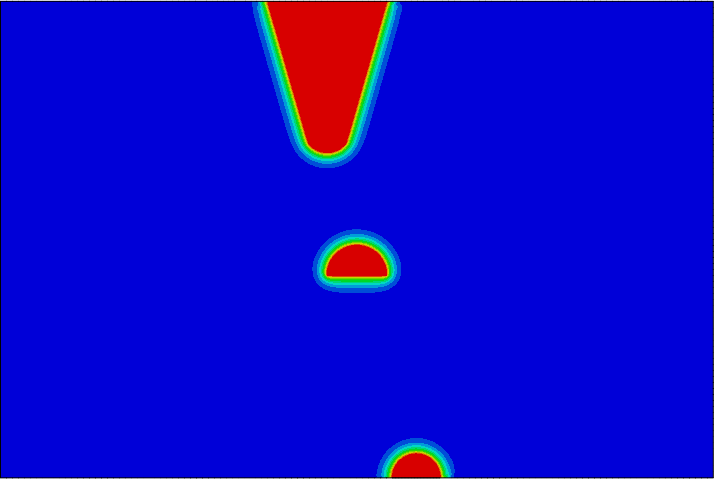}
                \caption{}
                \label{fig:Crack1}
        \end{subfigure}
        \begin{subfigure}[h]{0.6\textwidth}
                \centering
                \includegraphics[scale=0.4]{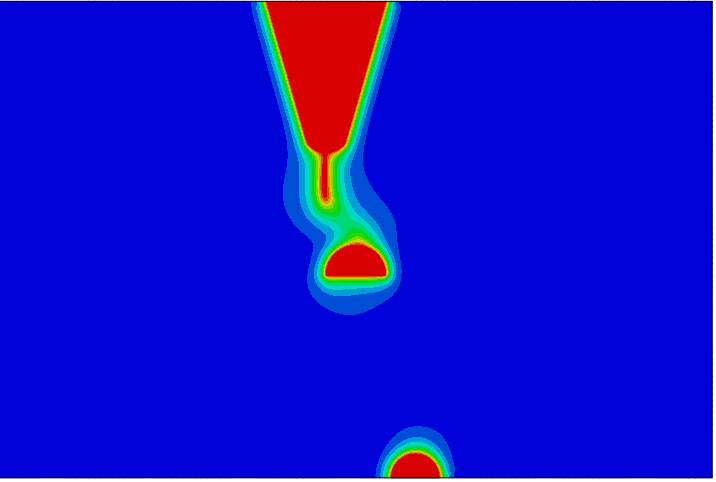}
                \caption{}
                \label{fig:Crack2}
        \end{subfigure}}
        
\makebox[\linewidth][c]{%
        \begin{subfigure}[h]{0.6\textwidth}
                \centering
                \includegraphics[scale=0.4]{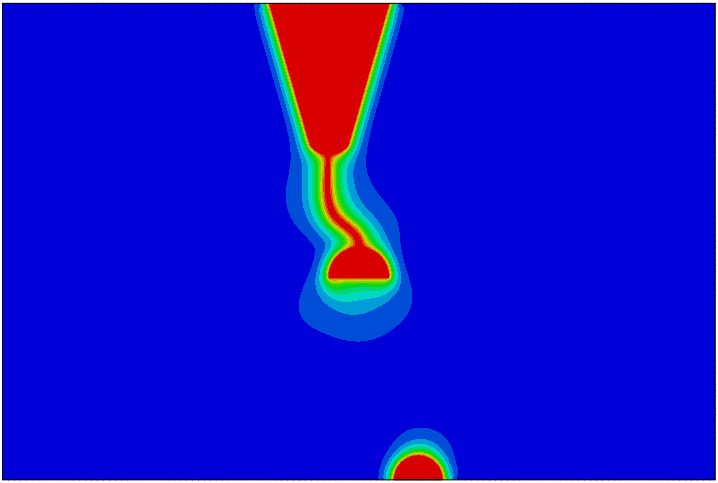}
                \caption{}
                \label{fig:Crack3}
        \end{subfigure}
        \begin{subfigure}[h]{0.6\textwidth}
                \centering
                \includegraphics[scale=0.4]{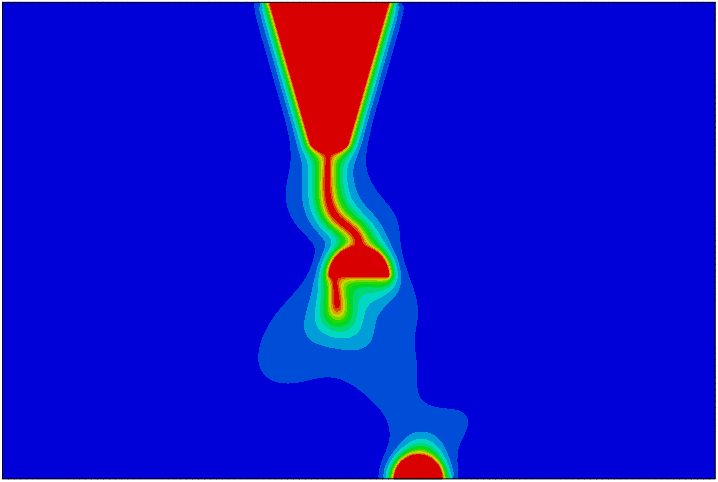}
                \caption{}
                \label{fig:Crack4}
        \end{subfigure}} 
        
\makebox[\linewidth][c]{%
        \begin{subfigure}[h]{0.6\textwidth}
                \centering
                \includegraphics[scale=0.4]{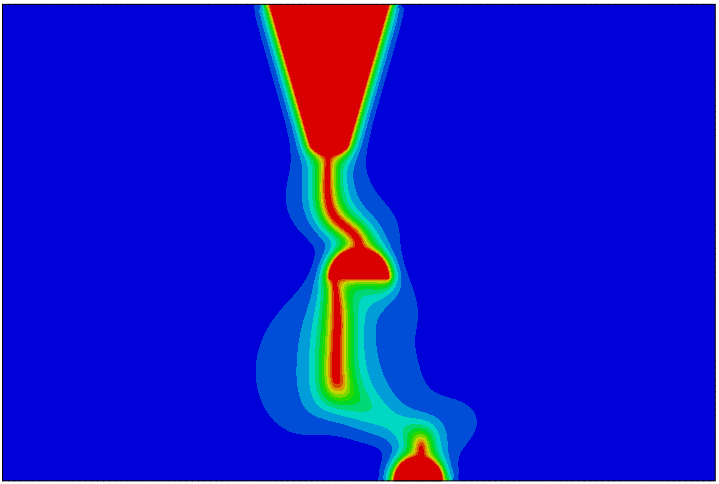}
                \caption{}
                \label{fig:Crack5}
        \end{subfigure}
        \begin{subfigure}[h]{0.6\textwidth}
                \centering
                \includegraphics[scale=0.4]{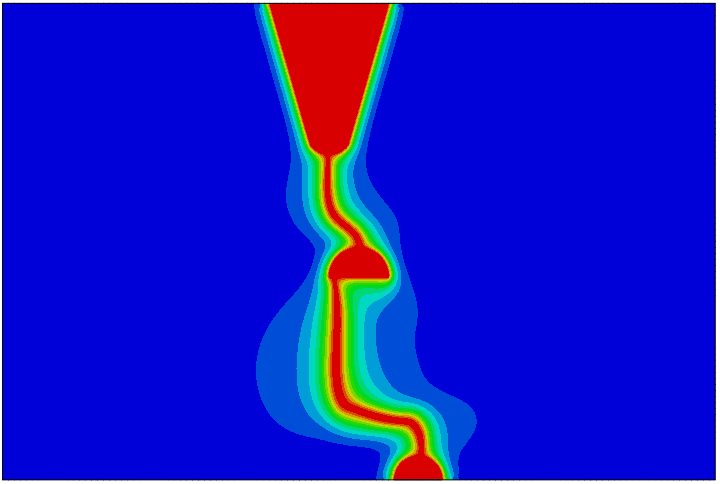}
                \caption{}
                \label{fig:Crack6}
        \end{subfigure}}        
       
        \caption{Rectangular plate with pre-existing corrosion defects. Fracture pattern at times (a) 1 s, (b) 5040 s, (c) 5721 s, (d) 6322 s, (e) 6948 s, and (f) 7560 s.}\label{fig:CrackCorrosion}
\end{figure}

The crack propagation contours are shown as a function of time in Fig. \ref{fig:CrackCorrosion}. Blue and red colors correspond to the completely intact and the fully broken state of the material, respectively. First, damage is restricted to the corrosion pits, which act as stress risers and attract the hydrogen (Fig. \ref{fig:Crack1}). Then, cracking initiates at the largest defect (Fig. \ref{fig:Crack2}) and eventually coalesces with the central hole (Fig. \ref{fig:Crack3}). After arresting at the center defect, the crack reappears at its left corner (Fig. \ref{fig:Crack4}) and continues to propagate. At the same time, cracking initiates at the smaller pit (Fig. \ref{fig:Crack5}) and both propagating cracks eventually coalescence (Fig. \ref{fig:Crack6}). Crack initiation takes place after 5040 s while the final coalescence occurs 2520 s later. This benchmark problem shows the robustness and capacity of the present phase field scheme to model complex cracking conditions that may arise in practical applications; crack deflection, arrest, restart and coalescence are successfully captured.

\section{Conclusions}
\label{Sec:Concluding remarks}

A phase field formulation for the long-standing problem of hydrogen assisted cracking has been presented. A coupled deformation-diffusion-phase field scheme has been developed in the framework of the finite element method. The model builds upon a hydrogen-dependent surface energy degradation law based on Density Functional Theory that adequately characterizes the sensitivity of the fracture energy to atomic hydrogen concentration. A wide variety of case studies have been addressed to benchmark the potential of the present phase field formulation for hydrogen embrittlement. First, the cracked plate benchmark is modeled to verify the numerical implementation and examine the influence of hydrogen in the load-displacement curve. Then, the failure stress in notched cylindrical bars is computed for a very wide range of environmental conditions, with finite element results showing a good agreement with the experiments. The model is also used to estimate exposure times for fracture in sea water under constant loads; numerical predictions enable to identify the cracking threshold without the need for time consuming laboratory tests. Numerical experiments in Compact Tension specimens are then conducted to compute internal hydrogen stress intensity thresholds as a function of hydrogen content. Finally, the capabilities of the model to capture complex crack propagation arising from defects intrinsic to corrosive environments are investigated. Computations show that phase field methods for fracture are a suitable tool to achieve the elusive goal of lifetime prediction of engineering components undergoing hydrogen assisted cracking.
 
\section{Acknowledgments}
\label{Sec:Acknowledgments}

E. Mart\'{\i}nez-Pa\~neda acknowledges financial support from the Ministry of Economy and Competitiveness of Spain through grant MAT2014-58738-C3 and the People Programme (Marie Curie Actions) of the European Union's Seventh Framework Programme (FP7/2007-2013) under REA grant agreement n$^{\circ}$ 609405 (COFUNDPostdocDTU). C.F. Niordson acknowledges financial support from the Danish Independent Research Fund (DFF-7017-00121).




\bibliographystyle{elsarticle-num}
\bibliography{library}

\end{document}